%% file: RP_March22_2020.tex
\documentclass[11pt]{article}
\usepackage{color}
\usepackage{amsmath,amsfonts,amssymb,amsthm,epsfig,epstopdf,titling,url,array}

\newcommand{\be}{\begin{equation}}
\newcommand{\ee}{\end{equation}}
\newcommand{\bes}{\begin{equation*}}
\newcommand{\ees}{\end{equation*}}
\newcommand{\beqn}{\begin{eqnarray}}
\newcommand{\eeqn}{\end{eqnarray}}
\newcommand{\beqns}{\begin{eqnarray*}}
\newcommand{\eeqns}{\end{eqnarray*}}
\newcommand{\lkr}{\left(}
\newcommand{\lkv}{\left[}
\newcommand{\rkv}{\right]}
\newcommand{\rkr}{\right)}
\newcommand{\lfi}{\left\{}
\newcommand{\rfi}{\right\}}

\newcommand{\fr}[1]{(\ref{#1})}

\newcommand{\del}{\delta}
\newcommand{\Del}{\Delta}

\newcommand{\al}{\alpha}

\newcommand{\eps}{\epsilon}

\newcommand{\ga}{\gamma}
\newcommand{\te}{\theta}
\newcommand{\om}{\omega}
\newcommand{\lam}{\lambda}
\newcommand{\Up}{\Upsilon}

\newcommand{\sig}{\sigma}

\newcommand{\Om}{\Omega}

\newcommand{\EE}{\ensuremath{{\mathbb E}}}

\newcommand{\II}{\ensuremath{{\mathbb I}}}

\newcommand{\PP}{\ensuremath{{\mathbb P}}}

\newcommand{\RR}{\ensuremath{{\mathbb R}}}

\newcommand{\vect}{\mbox{vec}}
\newcommand{\Pen}{\mbox{Pen}}
\newcommand{\tPen}{\overline{\rm Pen}}
\newcommand{\Span}{\mbox{Span}}

\newcommand{\diag}{\mbox{diag}}

\newcommand{\etal}{{\it et  al. }}
\newcommand{\std}{\mbox{std}}

\newcommand{\Tr}{{\rm Tr}}

\newtheorem{theorem}{Theorem}
\newtheorem{lemma}{Lemma}
\newtheorem{corollary}{Corollary}

\newtheorem{remark}{Remark}

\newcommand{\ba}{\mathbf{a}}

\newcommand{\bof}{\mathbf{f}}
\newcommand{\bg}{\mathbf{g}}

\newcommand{\bh}{\mathbf{h}}

\newcommand{\bu}{\mathbf{u}}

\newcommand{\bw}{\mathbf{w}}

\newcommand{\by}{\mathbf{y}}

\newcommand{\bA}{\mathbf{A}}

\newcommand{\bD}{\mathbf{D}}
\newcommand{\bE}{\mathbf{E}}
\newcommand{\bF}{\mathbf{F}}
\newcommand{\bG}{\mathbf{G}}

\newcommand{\bI}{\mathbf{I}}

\newcommand{\bQ}{\mathbf{Q}}
\newcommand{\bS}{\mathbf{S}}

\newcommand{\bW}{\mathbf{W}}
\newcommand{\bX}{\mathbf{X}}
\newcommand{\bY}{\mathbf{Y}}
\newcommand{\bZ}{\mathbf{Z}}
\newcommand{\bK}{\mathbf{K}}

\newcommand{\bzero}{\mathbf{0}}

\newcommand{\boeta}{\mbox{\mathversion{bold}$\eta$}}

\newcommand{\beps}{\mbox{\mathversion{bold}$\eps$}}

\newcommand{\bGa}{\mbox{\mathversion{bold}$\Gamma$}}
\newcommand{\bom}{\mbox{\mathversion{bold}$\om$}}

\newcommand{\bUp}{\mbox{\mathversion{bold}$\Up$}}

\newcommand{\bSig}{\mbox{\mathversion{bold}$\Sigma$}}
\newcommand{\bTe}{\mbox{\mathversion{bold}$\Theta$}}

\newcommand{\bPi}{\mbox{\mathversion{bold}$\Pi$}}

\newcommand{\hbZ}{\widehat{\bZ}}
\newcommand{\tbZ}{\tilde{\bZ}}
\newcommand{\hJ}{\hat{J}}
\newcommand{\tJ}{\tilde{J}}
\newcommand{\hK}{\hat{K}}
\newcommand{\tK}{\tilde{K}}
\newcommand{\hfm}{\hat{f}_m}
\newcommand{\hbG}{\widehat{\bG}}
\newcommand{\hbf}{\hat{\bof}}
\newcommand{\tbf}{\tilde{\bof}}

\newcommand{\calA}{{\mathcal{A}}}
\newcommand{\calB}{{\mathcal{B}}}

\newcommand{\calG}{{\mathcal G}}
\newcommand{\calH}{{\cal H}}

\newcommand{\calM}{{\mathcal M}}

\newcommand{\calS}{{\mathcal{S}}}

\newcommand{\calW}{{\cal W}}

\newcommand{\calZ}{{\cal{Z}}}

\newcommand{\lan}{\langle}
\newcommand{\ran}{\rangle}

\newcommand{\di}{\displaystyle}

\long\def\ignore#1{}

\title{\Large{\bf Is clustering advantageous in statistical ill-posed linear inverse problems}}
\author{
\large{ Rasika Rajapakshage and Marianna Pensky} \\ 
  \\  
Department of Mathematics, University of Central Florida } 

\date{}

\begin{document}
\maketitle

\begin{abstract}
In many statistical  linear inverse problems, one needs to recover classes of similar objects 
from their noisy images under an operator that does not have a bounded inverse. 
Problems of this kind appear in many areas of application. Routinely, in such problems  
clustering is carried out at a pre-processing step and then the inverse problem is solved for 
each of the cluster averages separately. As a result, the errors of the procedures are usually 
examined for the estimation  step only. The objective of this paper is to examine, both theoretically 
and via simulations,  the effect of clustering on the accuracy of the solutions of general ill-posed 
linear inverse problems. 
In particular, we assume that one observes 
$X_m = A f_m + \del \eps_m$, $m=1, \cdots, M$, where functions $f_m$ can be grouped into $K$ classes
and one needs to recover a vector function $\bof= (f_1,\cdots, f_M)^T$. 
We construct  an estimator  for $\bof$ as a solution of  a penalized optimization problem
which corresponds to the clustering before estimation setting.  
We derive an oracle inequality for its precision and confirm that the estimator is 
minimax optimal or nearly minimax optimal up to a  logarithmic factor of the number of observations. 
One of the advantages of our approach  is that we do not assume that the number of clusters is 
known in advance. Subsequently, we compare the accuracy of the above procedure with the precision of estimation 
without clustering, and clustering following the recovery of each of the unknown functions separately.

We conclude that clustering at the pre-processing step is beneficial when the problem is moderately ill-posed.
It should be applied with extreme care when the problem is severely ill-posed. 
\\


\end{abstract}

\section {Introduction }
\label{sec:introduction}
\setcounter{equation}{0}

In this paper, we consider   a set of general ill-posed linear inverse problems $A f_m = q_m$, $m=1, \cdots, M$, where $A$ is a  
bounded linear  operator that  does not have a bounded inverse and the right-hand sides $q_m$ are measured with error. 
In particular, we assume that some of the objects   $f_m$ and hence  $q_m$, are very similar to each other, 
so that they can be averaged and recovered together. As a result, one supposedly  obtains estimators of $f_j$ with smaller errors. 
The grouping is usually unknown (as well as the number of groups)  and is carried out at a pre-processing step by applying one of the 
standard clustering techniques with the number of clusters determined by trial and error.  
Subsequently, the objects in the same cluster are averaged and the errors of those aggregated curves are used as true errors in the analysis.

Problems of this kind appear in many areas of application such as astronomy (blurred images), econometrics 
(instrumental variables), medical imaging (tomography, dynamic contrast enhanced Computerized Tomography  and Magnetic Resonance Imaging), 
finance (model calibration of volatility) and many others where similar objects are measured and can be recovered together.  
Indeed, clustering has been applied for decades  to solve ill-posed inverse problems in pattern recognition 
\cite{bezdek}, astronomy \cite{starck}, astrophysics \cite{fraix}, pattern-based time series segmentation \cite{deng}, 
medical imaging \cite{com1}, elastography for computation of the unknown stiffness distribution \cite{arnold} 
and  for detecting early warning signs on stock market bubbles \cite{kurum}, to name a few.
While in some  other  settings  the main objective is finding group assignments, 
we are  considering only applications where clustering is used merely as a denoising technique.
 In those applications, routinely, clustering is carried out at the pre-processing step 
and then the inverse problems are solved for each of the cluster averages separately. As a result, the errors of 
the procedures are usually examined for the estimation  step only.
The objective of this paper is to examine, both theoretically and via simulations,  
the effect of clustering on the accuracy of the solutions of general ill-posed linear inverse problems.

There exists immense literature on the statistical inverse problems (see, e.g.,  \cite{abram}, \cite{abr_pen},   \cite{bissantz}, 
\cite{blanchard}, \cite{cohen}, \cite{donoho}, \cite{pen_lasso}  and monographs  \cite{alquier}, \cite{engl} and references therein, to name a few).
However, to  the best of our knowledge, the question about the effects of clustering in  statistical inverse problems has never been investigated.
Recently, as a part of a more general  theory, the effect of clustering on the precision of recovery in multiple regression problems 
has been studied in \cite{klopp}.  Klopp \etal \cite{klopp} concluded that, even under uncertainty,  clustering improves 
the estimation accuracy. The goal of this paper is to extend this study to the ill-posed linear inverse problems setting.

In particular, we consider the following problem. 
Let $A: \calH_1 \to \calH_2$ be a known linear operator where  $\calH_1$ and $\calH_2$ are Hilbert spaces  
with inner products $\lan \cdot, \cdot \ran_{\calH_1}$ and $\lan \cdot, \cdot \ran_{\calH_2}$, respectively. 
The objective is to recover functions $f_m  \in {\cal H}_1$   from 
\be\label{eq1}
X_m (x) = q_m (x)   +  \del\, \epsilon_m (x), \quad q_m = A f_m, \quad m=1, \cdots, M,
\ee
where  $\eps_m (x)$ are the independent   white noise processes and the goal is to recover the vector function $f = (f_1, \cdots, f_M)^T$. 
Assume that observations are  taken as  functionals of   $X_m$: for any $\psi \in {\cal H}_2$ one observes
\be \label{eq:observ}
\lan X_m, \psi \ran = \lan  A f_m, \psi \ran +  \del\,  \xi_m(\psi),  
\ee
where $\del$ is noise level and $\xi_m(\psi)$ are zero  mean Gaussian random variables with   
\be \label{eq:cov}
\EE [\xi_m(\psi_1)  \xi_l(\psi_2)] = \lfi
\begin{array}{ll}
\lan \psi_1, \psi_2 \ran_{\calH_2},   & m=l\\
0, & m \neq l 
\end{array} \right.
\ee

In what follows we consider the situation where, despite of $M$ being large, there are only $K$ types of functions $f_m(t)$.
In particular, we assume that there exists a collection of functions $h_1(t),...,h_K(t)$ such that  $f_m(t) = h_k(t)$ 
for any $m$ and some $k =  z(m)$. In other words, one can define a clustering function $z = z(m)$, $m=1, \ldots, M$,  with values in  
$\{1, \ldots, K\}$ such that $f_m = h_{z(m)}$. We denote the clustering matrix corresponding to the clustering function $z(m)$  by $\bZ$.
Note that $\bZ \in \{0,1\}^{M\times K}$ and $\bZ_{m,k} = 1$ if and only if $z(m)=k$, so that matrix $\bD^2 = \bZ^T \bZ$ is diagonal. 

If   the function $z(m)$ were known, one could improve precision of estimating $f_m$ by averaging the signals 
within  clusters and construct the estimators $\hat{h}_k$ of the common cluster means, thus  reducing the noise levels,  
and subsequently set  $\hat{f}_m = \hat{h}_{z(m)}$.
In reality, however, neither the true clustering matrix $\bZ_*$, nor the true number of classes $K_*$  are  available, so 
they also need to be estimated.

Note that the  objective is accurate estimation of functions $f_m$, $m=1, \cdots, M$, rather than recovery 
of  the clustering matrix $\bZ$. Moreover, although a true clustering matrix $\bZ_*$ always exists (if all 
functions $f_m$ are different, one can choose $K_* = M$ and $\bZ_* = \bI_M$), one is not interested in finding 
$\bZ_*$. Indeed,  one would rather  incur a small bias resulting from  replacement of $f_m$ by $h_k \approx f_m$    than  
obtain estimators with high variances, that are common in inverse problems where each function $f_m$ is estimated separately. 
On the other hand,    using the clustering procedure leads to  one more type of errors that are due to erroneously
pooling together estimators of functions $f_m$ that belong to different classes, i.e., the errors due to mistakes in clustering.

The goal of this paper is the study of the theoretical recovery limits for the unknown functions $f_m$, $m=1, \cdots, M$,
when one applies clustering, thus taking advantage of the fact that some of the functions $f_m$ are similar to each other,
or ignores this knowledge and proceeds with estimation without clustering.
In order to evaluate benefits of clustering, we formulate estimation with clustering problem as an optimization problem.
One of the advantages of our approach  is that we do not assume that the number of clusters is known in advance.
Instead, we elicit  the unknown number of clusters,   the clustering matrix and the estimators of the unknown functions as a solution 
of  a penalized optimization problem where a penalty is placed on the unknown number of clusters. 
For this reason, our analysis applies not only to an ``ideal'' (but usually impractical) situation when the number of 
clusters is known but to the realistic scenario when it is unknown.

In this paper we analyze the situation where clustering is done before estimation, at the pre-processing level, as it usually
happens in many applications. The optimization problem in the paper corresponds to this scenario (specifically, to the K-means clustering setting),
as well as our in-depth theoretical study which evaluates the precision of estimators with   
clustering and compares it to the estimation accuracy without clustering. 
In order to further assess  benefits of clustering, we implement  a numerical study and 
compare the estimators  where clustering was carried out at the pre-processing level 
(``Clustering before'') to the estimators where clustering was done post-estimation (``Clustering after'') 
and the estimators  without clustering (``No clustering'').
We conclude   that clustering at the pre-processing level improves estimation precision when the inverse 
problem is moderately ill-posed but brings no benefits (and can even increase estimation errors) if the 
problem is severely ill-posed.

The rest  of the paper is organized as follows. In Section \ref{sec:assump_est}, we introduce notations and assumptions and 
discuss   optimization problem  that delivers the estimator. Section \ref{sec:error} deals with 
quantification of estimation errors. In particular, Section \ref{sec:oracle} provides the oracle 
expression for the risk of an estimator obtained in Section \ref{sec:estimation}. Section \ref{sec:minimax_upper}
presents  upper bounds for the risk under the assumptions in Section \ref{sec:assump}. 
In order to ensure that the estimators in Section \ref{sec:estimation} are asymptotically optimal, 
in Section \ref{sec:minimax_lower} we derive   minimax  lower bounds for the risk. Finally,
Section \ref{sec:advantage} carries out theoretical comparison of estimation accuracy with and without clustering
in asymptotic setting. Section \ref{sec:simulations} performs a similar comparison via a simulation study 
for the case of finite-valued parameters. 
 Finally, Section \ref{sec:discussion} contains in-depth discussion and recommendations about application of the pre-clustering 
in the linear ill-posed problems.    
Section~\ref{sec:proofs} contains proofs of all statements in the paper.


\section{Assumptions and estimation }
\label{sec:assump_est}
\setcounter{equation}{0} 

\subsection{Notations}

Below, we shall use the following notations. We denote $[m] = \{1, \cdots,m\}$.
We denote vectors and matrices by bold letters. 
For any vector $\ba$, we denote its $l_2$-norm by $\|\ba\|$ and the $l_0$ norm, the number of non zero elements, by $\|\ba\|_0$ . 
For any matrix $\bA$, we denote its Frobenius norm by $\|\bA \|_F$, the operator norm by $\|\bA \|_{op}$ 
and the span of the column space of matrix $\bA$ by $\Span(\bA)$. 
We denote the Hamming distance between matrices $\bA_1$ and $\bA_2$, the number of nonzero elements in $\bA_1-\bA_2$,  by 
$\|\bA_1 - \bA_2\|_H$. 
We denote the $(k \times k)$ identity matrix by $\bI_k$ and drop subscript $k$ when there is no uncertainty about the dimension.
We denote the inner product and the corresponding norm in a Hilbert space $\calH$ by $\lan \cdot, \cdot \ran_{\calH}$ 
and $\|\cdot\|_{\calH}$, respectively, and drop subscript $\calH$ whenever there is no ambiguity.
For any set $S$, we denote cardinality of $S$ by $|S|$.
We denote the set of all clustering matrices for grouping $M$ objects into $K$ classes by $\calM(M,K)$.
We denote  $a_n \lesssim b_n$ if there exist $c  < \infty$ independent of $n$ such that $ a_n \leq c b_n$ and 
$a_n \gtrsim b_n$ if there exist $c  > 0$ independent of $n$ such that $ a_n \geq c b_n$. Also, 
$a_n \asymp b_n$ if simultaneously  $a_n \lesssim b_n$ and $a_n \gtrsim b_n$.
Finally, we use   $C$ as a generic absolute constant independent of $n$, $M$ and $K$, which can take different values in different places.


\subsection{Reduction to the matrix model}
\label{sec:seq_model}

Since observations are taken as linear functionals \eqref{eq:observ},  the problem can be reduced to the so-called sequence model.
For this purpose, the unknown functions are expanded over an orthonormal basis  $\phi_j$, $j=1,2, \cdots,$ of $\calH_1$
and the problem reduces to the recovery of the unknown coefficients of those functions. This is a common technique in the field of 
statistical inverse problems (see, e.g., Cavalier {\it et al.} (2002), Cavalier and Golubev (2006)  and Knapik {\it et al.} (2011)). 
The orthonormal basis is commonly taken to be the eigenbasis of the operator $A$. 
Since the eigenbasis is often unknown, in this paper, we consider a wider variety of basis functions. Specifically,  
we assume that operator $A$ allows a wavelet-vaguelette decomposition introduced by Donoho (1995). In particular, 
Donoho (1995) assumed that  there exists an orthonormal basis $\phi_j$, $j=1,2, \cdots,$ of $\calH_1$
and nearly orthogonal sets of functions $\psi_j, \eta_j \in \calH_2$, $j=1,2, \cdots$, such that for some constants $\nu_j>0$,
and some absolute constants $0< c_\psi, C_\psi, c_\eta, C_\eta < \infty$ independent of $j$, 
one has for any vector $\ba$:
\begin{align}
& A \phi_j = \nu_j^{-1} \eta_j,\quad A^* \psi_j = \nu_j^{-1} \phi_j; 
\quad  \langle \eta_{j_1},\psi_{j_2}\rangle_{\calH_2} = I(j_1=j_2); \label{Don1}\\
& c_\psi^2 \|\ba\|^2 \leq \|\sum_j a_j\psi_j\|^2  \leq C_\psi^2 \|\ba\|^2,\quad 
c_\eta^2 \|\ba\|^2 \leq  \|\sum_j a_j \eta_j\|^2 \leq C_\eta^2 \|\ba\|^2, \label{Don3}
\end{align} 
where $A^*: \calH_2 \to \calH_1$  is the linear operator conjugate to $A$ and $I(\ldots)$ is the indicator function.
The name was motivated by the fact that conditions \fr{Don1} and \fr{Don3} hold for a variety of linear operators  such as 
convolution, numerical differentiation or Radon transform when  $\{\phi_j\}$ is a wavelet basis (see also 
 Abramovich and Silverman  (1998)).  
Obviously, assumptions \eqref{Don1} and \eqref{Don3} are valid when $\{\phi_j\}$ is the eigenbasis of the operator $A$. 
Under  conditions \eqref{Don1} and \eqref{Don3}, any function $f$ can be recovered from its image $Af$ using  reproducing formula
\be \label{eq:repr1}
f = \sum_{j} \nu_j\, \lan Af, \psi_j \ran \phi_{j}
\ee
which is analogous to the reproducing formula for the eigenbasis case.

We expand functions $f_m \in \calH_1$ over the   basis $\phi_j$,   $j=1, \cdots, $   and denote the matrix of coefficients by $\bG$.
Denote $\lan A f_m, \psi_j \ran =   \bQ_{j,m}$, so that, by  \eqref{eq:repr1}, for  $j=1,2, \cdots,$  $m=1, \cdots, M$, one has
\be \label{main_relation}
\bG_{j,m} = \lan f_m, \phi_j \ran = \nu_j  \lan f_m, A^* \psi_j \ran = \nu_j \,  \lan A f_m, \psi_j \ran =  \nu_j\,  \bQ_{j,m}.
\ee
Consider matrix of observations $\bY$ and matrix of errors $\bE$ with respective  components $\bY_{j,m} = \lan X_m, \psi_j \ran$ and 
$\bE_{j,m} = \xi_m(\psi_j)$ where $\xi_m(\psi)$ is defined in \eqref{eq:cov}. 
Let $\bG_*$ and $\bQ_*$ be the true matrices of coefficients. Then, it follows from 
\fr{eq1}, \fr{eq:observ} and \fr{main_relation} that elements $\bY_{j,m}$ of column $m$ of matrix 
$\bY$ obey  the sequence model 
\be \label{main_seq}
\bY_{j,m} = \nu_j^{-1} (\bG_*)_{j,m} + \del \bE_{j,m}, \quad j=1,2, \cdots, \quad m=1, \cdots, M.
\ee 
Here,  $\EE(\bE_{j,m}) =0$ and, by \fr{eq:cov}, 
\be \label{er_covar}
\EE (\bE_{j_1,m_1} \bE_{j_2, m_2}) = \lfi
\begin{array}{ll}
0, & m_1 \neq m_2\\
\lan \psi_{j_1}, \psi_{j_2} \ran, & m_1 = m_2
\end{array} \right.
\ee
In order to make the model computationally convenient, we cut the sequence model at some index  $n$ where $n$ is large enough 
to make the error, which is due to this reduction, negligibly small. 
Then, $j=1,\ldots, n$,  and $\bG_*$, $\bQ_*$, $\bY$ and $\bE$ are 
$n \times M$ matrices, Also,  it follows from \eqref{main_seq} that 
\be \label{main_eq1}
\bUp  \bY  = \bG_* + \del \bUp  \bE, \quad \bUp = \diag(\nu_1,   \cdots, \nu_n).
\ee 
We shall discuss the choice of $n$ later in Section~\ref{sec:assump}.

Denote the matrix with elements $\bSig_{i,j} = \lan \psi_i, \psi_j \ran$  by $\bSig$ 
and observe that \eqref{er_covar} implies that
\be \label{eq:mom_matr_normal}  
\EE[(\bE\bE^T)]=M\, \bSig, \qquad \EE(\bE^T\bE)=n\, \bI_M.
\ee  
Hence, matrix $\bE$ has the matrix-variate normal distribution $\bE \sim N(0,\bSig \otimes \bI_M)$.
Observe that the first relation in formula \eqref{Don3}  implies that
\be \label{Sigma_norm}
\|\bSig\|_{op}    \leq C_\psi^2.
 \ee


\subsection{Assumptions}
\label{sec:assump}

Recall that functions $f_m$ belong to $K$ different groups, so that $f_m = h_k$ with $k = z(m)$ where 
$z = z(m)$ is a clustering function.  Denote the matrix of coefficients of functions $h_k$ in the basis $\phi_j$ by $\bTe$, so that
$\bTe_{j,k} = \lan h_k, \phi_j \ran$, $j=1, \cdots,n$, $k=1, \cdots, K$.

It is well known that recovery of an unknown function  from noisy observations  relies on the fact that it possesses some minimal 
level of smoothness. This  smoothness usually manifests as gradual decline  of coefficients of this function in some basis, 
so the coefficients  decrease as one uses more and more complex basis functions. 
For this reason, we assume that $h_k$ belong to a   ball:  $h_k \in \calS(r,\calA)$, $k=1, \ldots, K$, 
where 
\be \label{class_calB}
\calS(r,\calA) = \lfi h= \sum_j \te_j \phi_j:\ \sum_{j=1}^\infty |\te_j|^{2} j^{2  r} \leq \calA^{2} \rfi. 
\ee
If  $\phi_j$ is the Fourier  basis, then \eqref{class_calB} defines a well known Sobolev ball. 
Formula \eqref{class_calB}  implies that 
\be \label{coef_cond}
\sum_{j=1}^\infty |\bTe_{j,k}|^{2} j^{2 r} \leq  \calA^{2}, \quad k=1, \ldots, K. 
\ee
If $r \geq 1/2$, then one can set the cut-off value to $n \approx \del^{-2}$. Indeed, 
the error rate in the problem cannot be smaller than a parametric rate of $C \del^{2}$
and \eqref{coef_cond} implies that the approximation error with this  value of $n$ will not exceed 
\be \label{eq:tail_cond}
\sum_{j=n+1}^\infty |\bTe_{j,k}|^{2} \leq n^{-2r} \sum_{j=1}^\infty |\bTe_{j,k}|^{2} j^{2 r} 
\leq \calA^{2} n^{-2r} \leq \calA^{2} \del^2
\ee  
In addition, it is well known (\cite{tsybakov}) that, in the regression setting,  
the observational version of the white  noise model \eqref{eq1}
based on a sample of size $n$ leads to $\del = \sig/\sqrt{n}$ where $\sig$ is the standard deviation of the noise.

Furthermore, since operator $A$ does not have a bounded inverse, the values of $\nu_j$ in \eqref{Don1} are growing with $j$.
While one can consider various scenarios, the standard assumption  is 
that $\nu_j$ grow monotonically with $j$  (see, e.g., Alquier {\it et al.}  (2011)):
\be \label{nu_j_cond}
\aleph_1   j^{\ga} \exp\lkr  \alpha j^\beta  \rkr \leq |\nu_j| \leq 
\aleph_2   j^{\ga} \exp\lkr  \alpha j^\beta  \rkr 
\ee
for some absolute positive constants $\aleph_1$, $\aleph_2$ and nonnegative $\ga$, $\alpha$ and $\beta$
where $\beta =0$ and $\ga >0$ whenever $\alpha =0$. 
The problem \eqref{eq1} is called {\it moderately ill-posed} if $\alpha=0$ and {\it severely ill-posed} if $\alpha>0$.


\subsection{Clustering and estimation}
\label{sec:estimation}

In what follows, we denote the true quantities using the star  symbol, i.e., 
$K_*$ is the true number of clusters, $\bZ_*$ is the true clustering matrix, 
$\bG_*$, $\bQ_*$ and $\bTe_*$ are the true versions of matrices  
 $\bG$, $\bQ$ and $\bTe$ and so on. As it was indicated before, we choose $n = [\del^{-2}]$, 
the largest integer that is no greater than $\del^{-2}$.

If  $z: [M] \to [K]$ is the clustering function and 
$\bZ \in \{0,1\}^{M\times K}$  is a clustering matrix, then  $\bG_{i,j} = \bTe_{i,z(j)}$ for
$i=1, \ldots, n$,  $j=1, \ldots, M$. Therefore, if the clustering matrix $\bZ$ were known, then 
one would repeat columns of matrix $\bTe$ to obtain $\bG$ and average columns of $\bG$ to construct $\bTe$.
Specifically, $\bG  = \bTe \bZ^T$ and $\bTe =   \bG \bZ \bD^{-2}$, 
where matrix $\bD^2 = \bZ^T \bZ$ is diagonal.

Denote by $\bPi_{\bZ,K}$ and $\bPi_{\bZ,K}^\bot$ the projection matrices on the column space of matrix $\bZ$ 
and    on the orthogonal subspace, respectively: 
\be  \label{eq:proj_def}
 \bPi_{\bZ,K} = \bZ (\bZ^T \bZ)^{-1}  \bZ^T, \quad \bPi_{\bZ,K}^\bot = \bI_M - \bPi_{\bZ,K}.
\ee 
Here, we use index $K$ to indicate that not only the clustering matrix $\bZ$ but also the number of clusters $K$ is unknown.
The projection matrix $\bPi_{\bZ,K}$ is such, that for any matrix $\bG \in \RR^{n \times M}$,   $\bG \bPi_{\bZ,K}$  
replaces each column of $\bG_j$ of $\bG$ by its average over all columns in cluster $z(j)$.
Then, matrix $\bG_*$ is such that  $\bG_*   = \bG_* \bPi_{\bZ_*,K_*}$ and, due to \eqref{main_eq1}, if $\bZ_*$ were known, 
it would seem  to be reasonable to estimate $\bG_*$  by $\bUp \bY \bPi_{\bZ_*,K_*}$. 
It is well known, however, that this estimator is inadmissible  and one needs to 
shrink or threshold elements of matrix  $\bUp \bY \bPi_{\bZ_*,K_*}$ to achieve an optimal 
bias-variance balance (\cite{mallat}, Section 11.2).

Observe that, since for the ill-posed inverse problems, the values of $\nu_j$ are growing with $j$ 
due to equation \eqref{nu_j_cond}, the elements $\bG_{j,i} = \bTe_{j,z(i)}$ of matrix  $G$ are harder and harder  to recover 
as $j$ is growing. On the other hand, condition \eqref{coef_cond} means that  coefficients $\bTe_{j,k}$ decrease rapidly 
as $j$ increases, and hence, for large $n$, one does not need to keep all $n$ coefficients 
for an accurate estimation of functions  $h_k$  (and therefore $f_m$). On the contrary, this will yield an estimator with a huge variance. 
For this reason, due to the fact that  conditions \eqref{coef_cond} apply to all $k = 1, \cdots, K$
simultaneously, we need to choose a set $J  \subseteq \{1, \ldots, n\}$  and set $\bTe_{jk}=0$ if $j \not\in J$. 
Then,  one has $\bG_{j,m} = 0$ if $j \in J^c$ where the set $J^c$ is complementary to $J$.
In order to express the latter in a matrix form, we introduce matrix 
\be \label{eq:bWJ}
\bW_J = \diag(\bw_1,...,\bw_n) \quad \mbox{with}\quad \bw_j = \II(j \in J), 
\ee
and observe that, for any matrix $\bG$, condition $(\bI_n - \bW_J)\bG=\bzero$ ensures that   $\bG_{j,m} = 0$, $j \in J^c$.

Consider integer  $K \in [M]$,  set  $\calM(M,K)$  of clustering matrices that cluster  $M$ nodes into $K$ groups
and set  $J \subseteq \{1, \ldots, n\}$.
Then, the objective is  to find matrices $\bG$ and $\bZ \in \calM(M,K)$, a set $J$ and an integer $K$:    
\begin{align}  \label{no_pen_opt0} 
  (\hat{\bZ}, \widehat{\bG},\hat{J}, \hat{K}) & \in \   \underset{\bZ,\bG,J,K}{\operatorname{argmin}}\,      
  \lfi  \|\bG - \bUp \bY \bPi_{\bZ,K} \|_F ^2 + \|\bUp \bY \bPi_{\bZ,K}^\bot \|_F ^2 \rfi   \\
&  \mbox{subject to} \  (\bI_n - \bW_J)\bG=\bzero, \nonumber
\end{align}
where $\bPi_{\bZ,K}^\bot$ is defined in \eqref{eq:proj_def}.
The second term in \eqref{no_pen_opt0} corresponds to the error of the $K$-means clustering of the matrix $\bUp \bY$
while the first term quantifies the difference between the clustered version of data matrix $\bUp \bY \bPi_{\bZ,K}$ 
and the matrix $\bG$.

Since $\|\bUp \bY \bPi_{\bZ,K} \|_F ^2 + \|\bUp \bY \bPi_{\bZ,K}^\bot \|_F ^2 = \|\bUp \bY   \|_F ^2$  
is independent of $\bG$ and $\bZ$,   the problem can be re-written in an equivalent form as 
\be \label{no_pen_opt}
 (\hat{\bZ}, \widehat{\bG},\hat{J}, \hat{K}) \in \   \underset{\bZ,\bG,J,K}{\operatorname{argmin}}\, \left\{
 \|\bG\|_F^2 -2\Tr (\bY^T \bUp \bG \bPi_{\bZ,K}) \rfi \ \mbox{subject to} 
\  (\bI_n - \bW_J)\bG=\bzero. 
\ee 
Note though that optimization problem    \fr{no_pen_opt} has a trivial  solution: 
$K=M$, $J = [n]$, $\bZ = \bI_M$ and $\bG = \bUp \bY$.

In order to avoid this, we   put a penalty on the value of $K$ and the set $J$, and 
find   $\bZ,\bG,J$ and $K$  as a solution of the following optimization problem:
\begin{align}\label{opt_prob1}
 (\hat{\bZ}, \widehat{\bG},\hat{J}, \hat{K}) \in \  & \underset{\bZ,\bG,J,K}{\operatorname{argmin}}\, \left\{    
 \|\bG\|_F^2 -2\Tr (\bY^T \bUp \bG \bPi_{\bZ,K} )+\Pen(J,K)\right\}\\
& \mbox{subject to}\  \bZ \in \calM(M,K),  (\bI_n- \bW_J)\bG=\bzero,  J \subseteq [n], K \in [M]. 
\nonumber
\end{align}
Optimization procedure \eqref{opt_prob1} leads to group thresholding of the rows of matrix $\bUp \bY \bPi_{\bZ,\bK}$
due to the condition $(\bI_n- \bW_J)\bG=\bzero$.  Indeed, if $\hat{\bZ}, \hat{J}$ and $\hat{K}$ were known, then it follows 
from \eqref{no_pen_opt0} that $\widehat{\bG}$ would be given by  
\be \label{eq:opt_bG}
\widehat{\bG} = \bW_{ \hat{J}} \bUp \bY \bPi_{\hat{\bZ},\hat{K}}
\ee
and problem \eqref{opt_prob1} can be presented as 
\begin{align} \label{opt_prob2}
 (\hat{\bZ},\hat{J}, \hat{K}) \in  \ &  \underset{\bZ,J,K}{\operatorname{argmin}}\, \left\{    
 \|(\bI- \bW_J)\bUp \bY \bPi_{\bZ,K} \|_F ^2  + \|\bUp \bY \bPi_{\bZ,K}^\bot \|_F ^2  +\Pen(J,K)\right\} \\
& \mbox{subject to}\  \bZ \in \calM(M,K),  
J \subseteq [n], K \in [M]. 
\nonumber
\end{align}
Note that the objective function  in \eqref{opt_prob2} is a sum of two  components: the first one 
is responsible for the best fitting of the matrix $\bUp \bY \bPi_{\bZ,K}$ when some of its rows are set to zero 
while the second one corresponds to the error of the $K$-means clustering of columns of matrix  $\bUp \bY$.
The solution of the optimization problem relies on the $K$-means algorithm that is NP-hard but, however, is known to 
provide very accurate results as long as initialization point is not too far from the true solution.

In practice, we shall solve optimization problem \eqref{opt_prob2} separately for each $K \in [M]$  
and then choose the value of $K$ that delivers the smallest value in \eqref{opt_prob2}.
We estimate the matrix of coefficients $\bG$ by  $\hbG$ defined in \eqref{eq:opt_bG}.
After coefficients $\hbG$ are obtained, we estimate $f_m$, $m=1, \ldots, M$, by 
\be \label{eq:hat_fm}
\hfm = \sum_{j \in J}  \widehat{\bG}_{j,m} \phi_j, \quad m=1, \cdots, M.
\ee

The penalty in  \eqref{opt_prob1} and \eqref{opt_prob2} should be chosen to exceed the random errors level
with high probability. If the number of clusters $K$, the set $J$ and the clustering matrix $\bZ$ were known, 
then the penalty would be of the order of the variance term $\di K  \sum_{j \in J} \nu_j ^2$. However, since
$K$, $J$ and $\bZ$ are unknown, we need to account for the uncertainty in estimation of those parameters by applying 
a union bound  and, hence, adding the terms that are proportional to the log-cardinality of the sets of those parameters. 
Since one has $n$ choices for $K$, $K^M$ possible clustering arrangements and approximately $\exp\lfi|J|\ln(ne/|J|)\rfi$
sets $J$ of cardinality $|J|$ for every $|J| = 1, ...,n$, we need to add a term proportional to
$\di (\max_{j \in J} \nu_j ^2)    \lkv M \ln K + |J|\ln(ne/|J|) + \ln (M n) \rkv.$ 
Finally, we need to choose a  constant $\tau$ and add a term proportional to $\di  \max_{j \in J} \nu_j ^2 \,  \ln (\del^{-\tau})$ to ensure that 
the upper bound holds with probability at least $1 - 2 \del^{\tau}$.
By carefully evaluating the upper bounds for each of the components of the error, we derive the penalty
\begin{align}
\begin{small}  
\label{eq:penalty} 
\Pen(J,K) =  2C_\psi^2 \del^2 \lkv 26 K  \sum_{j \in J} \nu_j ^2 +
39  (\max_{j \in J} \nu_j ^2)    \lfi M \ln K +|J|\ln \lkr \frac{ne}{|J|}\rkr  +\ln \lkr \frac{M n}{\del^\tau}\rkr   \rfi \rkv  
\end{small} 
\end{align}
where $n = [\del^{-2}]$, $C_\psi$ is defined in \eqref{Sigma_norm} and the choice of $\tau$ ensures that 
the upper bound for the error will hold with probability at least $1 - 2 \del^{\tau}$. 
Hence, in any real life setting, the constant $\tau$ should be such that 
this probability is large enough.

Penalty \eqref{eq:penalty} consists of  four  terms.
The first term, $26 K  \sum_{j \in J} \nu_j ^2$ represents the error of estimating $|J|$ 
coefficients for each of the distinct functions $h_k$, $k=1, \ldots, K$.
The second and the third terms account  for the difficulty of clustering $M$ functions into $K$ classes 
and choosing a set $J \subset \{1, \ldots, n \}$. The   last term   is of the smaller asymptotic order, it 
offsets the error of the choice of $K$ and also ensures that the oracle inequality holds 
with the probability at least $1 - 2 \del^{\tau}$. 
Observe that since the data is weighted by the diagonal matrix $\bUp$
in \eqref{main_eq1}, the last three terms are weighted by $\max_{j \in J} \nu_j ^2$.

The penalty \eqref{eq:penalty}  corresponds to the general model selection that does not rely on assumptions  
\eqref{class_calB} and \eqref{nu_j_cond}. If those conditions hold, the elements $(\bG_*)_{j,m}$ are decreasing with $j$  
for every $m$, while the values of $\nu_j$ are increasing. Therefore, one should choose 
a set $J$ of the form $J = \lfi 1, ..., L \rfi$ for some $L \leq n$. Since the cardinality of the set of possible $L$'s 
is just $n$, this would lead to replacement of the term $|J|\ln \lkr  ne/|J| \rkr$ in the penalty by merely $\ln n$
leading to 
\begin{align}
\begin{small}  
\label{eq:penalty_new} 
\tPen(L,K) =  2C_\psi^2 \del^2 \lkv 26 K  \sum_{j=1}^L \nu_j ^2 +
39    \nu_L ^2\,     \lfi M \ln K   +\ln \lkr \frac{M n}{\del^\tau}\rkr   \rfi \rkv  
\end{small} 
\end{align}
%
 

\begin{remark} \label{rem1}   
{\bf (Unknown noise level). \  }
{\rm 
The value of $\del$ in \eqref{eq:penalty} and \eqref{eq:penalty_new} is usually unknown but can be easily 
obtained from data. Indeed,   one can apply a wavelet transform to the original data matrix $\bY$,
and then recover $\del$ as the median of the absolute value of the wavelet coefficients  at the highest resolution level 
divided by 0.6745  (see, e.g., Mallat (2009), Section 11.3). In fact, in our simulations, we treated $\del$ as an unknown quantity and 
estimated it by this procedure.  
}
\end{remark}


\begin{remark} \label{rem2}   
{\bf (Different smoothness for different clusters).\  }
{\rm 
One can consider a more general case where functions from 
different clusters have different smoothness levels. 
In this case, each function $h_k$ has a corresponding set  of
nonzero coefficients $J_k$, $k=1, \ldots, K$, which may be of the form $\{1, \ldots, L_k\}$. Consequently,
the    terms $\di  K  \sum_{j \in J} \nu_j ^2$ and $\di  K  \sum_{j=1}^L \nu_j ^2$ in the penalties  
 \eqref{eq:penalty} and \eqref{eq:penalty_new} should be replaced by, respectively, 
\bes
\sum_{k=1}^K\, \sum_{j \in J_k} \nu_j ^2\ \quad \mbox{and} \quad \sum_{k=1}^K\, \sum_{j=1}^{L_k} \nu_j ^2.
\ees 
Theoretical results for this case are a matter of a future investigation.
}
\end{remark}


\section{Estimation error }
\label{sec:error}
\setcounter{equation}{0}


\subsection{The oracle inequality}
\label{sec:oracle}

The average error of estimating  $f_m$ by $\hfm$,  $m=1,\ldots, M,$ is  given by
\be \label{eq:fun_est_error}
R (\bof, \hbf) = M^{-1} \, \sum_{m=1}^M \|\hfm - f_m\|^2,
\ee
where $\bof$ and $\hbf$ are column vector  with functional components $f_m$ and $\hfm$,  
$m=1,\ldots, M,$ respectively. 
Due to the inequality  \eqref{eq:tail_cond}, the errors of approximation of functions $f_m$ by 
the $n$-term expansions over   $\phi_j$, $j=1,\ldots, n$, are much smaller 
than the errors due to estimation or thresholding of the first $n$ coefficients of these expansions.
Therefore, the main portion of the error is due to $M^{-1}  \|\widehat{\bG}-\bG_*\|_F ^2$.
The  following statement  places an upper bound on $\|\widehat{\bG}-\bG_*\|_F ^2$.

\begin{theorem} \label{th:upper_bound}
Let $(\hat{\bZ}, \widehat{\bG},\hat{J}, \hat{K})$ be a solution of optimization problem \eqref{opt_prob1}
with the penalty $\Pen(J,K)$  given by expression \eqref{eq:penalty}. 
Then, there exists a set $\Omega = \Omega(\tau)$ 
with   $\PP (\Om) \geq 1 - 2\del^{\tau}$  such that for every $\om \in \Om$ one has 
\begin{align} \label{eq:upper_bound}
\|\widehat{\bG}-\bG_*\|_F ^2  & \leq \min_{\bZ,J,K} \lfi 3\, \|\bW_J \bG_* \bPi_{\bZ,K}-\bG_*\|_F ^2 + 
4\, \Pen(J,K) \rfi
\end{align}
Moreover, if assumptions \eqref{class_calB} and \eqref{nu_j_cond} hold and  
$(\hat{\bZ}, \widehat{\bG},\hat{L}, \hat{K})$ is a solution of optimization problem \eqref{opt_prob1}
with  $J = \{ 1, ..., L \}$  and the penalty $\Pen(J,K)$ replaced with $\tPen(L,K)$ 
defined in  \eqref{eq:penalty_new}, then, for   $\om \in \Om$ 
\begin{align} \label{eq:upper_bound_new}
\|\widehat{\bG}-\bG_*\|_F ^2  & \leq \min_{\bZ,L,K} \lfi 3\, \|\bW_J \bG_* \bPi_{\bZ,K}-\bG_*\|_F ^2 + 
4\, \tPen(L,K) \rfi
\end{align}
\end{theorem}

Theorem~\ref{th:upper_bound}  provides an oracle inequality for  $\|\widehat{\bG}-\bG_*\|_F ^2$.
The first term in expression \eqref{eq:upper_bound} is the bias term that quantifies the error of approximation 
of matrix $\bG_*$ when its columns are averaged over $K$ clusters using matrix $\bZ$ and one keeps only  
terms with $j \in J$ in the approximations of each of the $K$ cluster means. This term is decreasing 
when $K$ and $|J|$ are increasing. The second term, $\Pen(J,K)$,  is the variance term that 
represents the error of estimation for the particular choices of $\bZ$, $J$ and $K$. This term 
grows when $K$ and $|J|$ are increasing. The error is provided by the best possible bias-variance balance in \eqref{eq:upper_bound}.

Since the right hand side in \eqref{eq:upper_bound} is minimized over $\bZ$ and $K$,    
if some of the functions $h_k$, $k=1\cdots, K$, are similar but not exactly identical 
to each other, it may be advantageous to place those functions in the same cluster, 
hence, reducing the variance component of the error. Our methodology will automatically 
take advantage of this opportunity.  Note that the error bounds in \eqref{eq:upper_bound} 
are non-asymptotic and are valid for any true matrix $\bG_*$ and any relationship between $K$, $M$ and $\delta$.

While those results are very valuable, they do not allow to quantify the effect of clustering on estimation errors 
when $\del$ is small and $M$ is large, so that $\del \to 0$, $M \to \infty$, and possibly $K \to \infty$. 
In the next section we shall investigate this issue under assumptions of Section~\ref{sec:assump}.


\subsection{The   upper bounds for the estimation error}
\label{sec:minimax_upper}

In order to study particular scenarios, in  what follows, we assume that $\nu_j$ satisfies condition \eqref{nu_j_cond}.
Assume that  $h_k \in \calS(r,\calA)$, $k=1, \ldots, K_*$,   where $\calS(r,\calA)$  is defined in \eqref{class_calB}.
Denote by $\bh$ the functional column vector with components $h_k$, $k=1, \ldots, K_*$.
Consider the maximum risk of our estimator $\hbf$ over all $h_k \in \calS(r,\calA)$, $k=1, \ldots, K_*$, and all true clustering matrices 
$\bZ_* \in \calM(M,K_*)$
\begin{align} \label{max_risk}
& R (\hbf,\calS(r,\calA),M, K_*)  = \max_{\bof,\bZ_*}\,   R (\bof, \hbf)   \quad   \mbox{subject to}\quad \\
& \bof = \bZ_*\, \bh, \ 
h_k \in \calS(r,\calA), \ k=1, \ldots, K_*, \ \bZ_* \in \calM(M,K_*), \nonumber
\end{align}
where $\calS(r,\calA)$ is defined in \eqref{class_calB} and $\calM(M,K_*)$ is the set of all 
clustering matrices that place $M$ objects into $K_*$ classes. 

In what follows, we assume that both $n$ and $M$ are growing simultaneously, that is, 
$\ln M \asymp \ln (n)$.  Note that this is a mild condition since it is satisfied when $M$ is 
growing at a rate of any positive power of $n$ or visa versa. Hence, due to $n \approx \del^{-2}$, we obtain
\be \label{eq:M_n_rel}
\ln (\del^{-1}) \asymp \ln n \asymp \ln M \asymp \ln (Mn).
\ee
Observe that  the first relation  follows from the definition of $n$  while the third one is the direct consequence of the second. 
Note also that the second assumption is both very mild and very natural. Since $\ln x$ grows very slowly with $x$,
in practical terms, it merely states that both $M$ and $\del^{-1}$ tend to infinity. The main consequence of the assumption 
\eqref{eq:M_n_rel}  is that the terms $\ln (\del^{-1})$,  $\ln n$ and $\ln M$ become interchangeable up to a constant.

Then, application of the oracle inequality  \eqref{eq:upper_bound} with $|J| =L$ and $K = K_*$ provides the following 
upper bounds for the error.

\begin{theorem} \label{th:upper_bound_L}
Let assumption \eqref{eq:M_n_rel} hold and    $\nu_j$, $j=1, \cdots, n$, satisfy condition \eqref{nu_j_cond} with  $r   \geq 1/2$.
  Let $(\hat{\bZ}, \widehat{\bG},\hat{L}, \hat{K})$ 
be a solution of optimization problem \eqref{opt_prob1} with the penalty given by expression \eqref{eq:penalty}.
Then, with probability at least $ 1 - 2\del^{\tau}$, one has 
\bes
R (\hbf,\calS(r,\calA),M,K_*) \leq C\,  R (M,K_*,\del),
\ees 
where the constant $C$ depends on $\alpha, \beta, \ga,   r, \tau$ and $\calA$ only and  
\be  \label{eq:RMKn_0}
R (M,K_*,\del) = 
\lkr  \del^2\, \ln K_* \rkr^{\frac{2r}{2r +2\ga}} + 
  \lkr {\del ^2\, M^{-1}  K_*}\rkr^{\frac{2r}{2r+2\ga +1}},     
\ee
if $ \alpha=\beta =0$, and 
\be  \label{eq:RMKn_1}
R (M,K_*,\del) = 
 \lkv \ln \lkr \frac{1}{ \del ^2   \ln K_*} \rkr \rkv ^{-\frac{2r}{\beta}} +  
\lkv \ln \lkr \frac{M }{ \del ^2    K_*} \rkr \rkv ^{-\frac{2r}{\beta}},  
\ee
if $\alpha>0, \beta > 0$.
\end{theorem}

\noindent
The expressions in \eqref{eq:RMKn_0} and \eqref{eq:RMKn_1} are well defined if 
$K_* \geq 2$. If $K_* =1$, then $\ln K_*=0$ and the first terms in \eqref{eq:RMKn_0}
and \eqref{eq:RMKn_1} are just equal to zero.


\subsection{The minimax lower bounds for the risk}
\label{sec:minimax_lower}

In order to show that the estimator developed in this paper is asymptotically near-optimal, 
below we derive minimax lower bounds for the risk over all $h_k \in \calS(r,\calA)$, $k=1, \ldots, K_*$, and all clustering matrices 
$\bZ_* \in \calM(M,K_*)$. For this purpose, we define the minimax risk as 
\be \label{Eq:minimax_risk}
R_{\min} (\calS(r,\calA),M, K_*) = \min_{\tbf} R (\tbf,\calS(r,\calA),M, K_*) 
\ee 
where $\tbf$ is any estimator of $\bof$ on the basis of matrix of observations $\bY$.

\begin{theorem} \label{th:lower_bound_L}
Let  $\nu_j$, $j=1, \cdots,$ satisfy condition \eqref{nu_j_cond}  and $r \geq 1/2$.  
Then, with probability at least $0.1$, 
one has 
\be \label{eq:low_bou_main}
R_{\min} (\calS(r,\calA),M, K_*) \geq C R_{\min} (M,K_*,\del)
\ee 
where the constant $C$ depends on $\alpha, \beta, \ga,   r$ and $\calA$ only and  
\be  \label{eq:RMKn_min_0}
R_{\min} (M,K_*,\del)  = 
\max \lfi \lkr {\del ^2\, \ln K_*}\rkr^{\frac{2r}{2r+2\ga}},\ 
\lkr  \del ^2\, M^{-1}\, K_* \rkr^{\frac{2r}{2r +2\ga +1}}
 \rfi, 
\ee
if $ \alpha=\beta =0$, and 
\be  \label{eq:RMKn_min_1}
R_{\min} (M,K_*,\del)  = 
\max \lfi \lkv \ln \lkr \frac{1}{ \del ^2 \ln  K_*} \rkr \rkv ^{-\frac{2r}{\beta}},\ 
\lkv \ln \lkr \frac{M }{ \del ^2    K_*} \rkr \rkv ^{-\frac{2r}{\beta}} \rfi,
\ee
if $\alpha>0, \beta > 0$.
\end{theorem}

\noindent
Observe that expressions for the upper and the lower bounds 
of the risk  \eqref{eq:RMKn_0} and \eqref{eq:RMKn_min_0} in the case of $\al = \beta = 0$, and 
\eqref{eq:RMKn_1} and \eqref{eq:RMKn_min_1} in the case of $\al >0, \beta > 0$
are identical, so our estimators are asymptotically optimal.


\subsection{The advantage of clustering}
\label{sec:advantage}

Theorems \ref{th:upper_bound_L} and \ref{th:lower_bound_L}  allow  to answer the question whether 
clustering in linear ill-posed inverse problems improves the estimation accuracy as $M \to \infty$ and $\del \to 0$.
Indeed,   solving problem \eqref{eq1} for each $m=1, \cdots, M$ separately 
is equivalent to choosing $K=M=1$ in the penalty. In this case,  one obtains the following corollary.

\begin{corollary} \label{cor:comparison}
If each of the inverse problems is solved separately, where the penalty is of the form \eqref{eq:penalty} 
with  $K=M=1$ and  $J = \lfi 1, \cdots, L \rfi$, then, with probability at least $ 1 - 2\, \del^{\tau}$, 
the average estimation error $\tilde{R} (\del)$ defined in \eqref{eq:fun_est_error} is bounded by
\be  \label{eq:noclust_upper_bound}
\tilde{R} (\del) \asymp \lfi 
\begin{array}{ll}
 \lkv  \del^2\rkv^{\frac{2r}{2\ga +2r +1}},
& \mbox{if}\ \alpha=\beta =0,\\
 \lkv \ln (\del^{-1}) \rkv^{-\frac{2r}{\beta}},
& \mbox{if}\ \alpha>0, \beta >0. 
\end{array} \right. 
\ee
If $r \geq 1/2$ and assumption \eqref{eq:M_n_rel} holds, then 
for $\del \to 0$, $M \to \infty$, one has  
\be  \label{eq:comparison}
\frac{R  (M,K_*,\del)} {\tilde{R} (\del)}
 \asymp \lfi 
\begin{array}{ll}
1 
& \mbox{if}\ \alpha>0, \beta >0,   \\
M^{-\frac{2r}{2\ga +2r +1}},
& \mbox{if}\ \alpha=\beta =0, K_*= 1\\
 \lkr\frac{K_*}{M}\rkr^{\frac{2r}{2\ga +2r +1}} + \lkr\del^2 \rkr^{\frac{2r}{(2\ga +2r +1)(2r + 2\ga)}} \ln(K_*),  
%
& \mbox{if}\ \alpha=\beta =0, K_*\geq 2.  
\end{array} \right. 
\ee 
Therefore, when $\del \to 0, \ M \to \infty$, clustering is asymptotically advantageous if $\alpha=\beta =0$.
\end{corollary}


\section {Simulations}
\label{sec:simulations}
\setcounter{equation}{0}


In order to study finite sample properties of the proposed estimation procedure, we carried out a   numerical study. 
In particular, we considered  a periodic convolution equation  $q = Ah = h * g$ with a kernel $g$ 
that transforms into a product in the Fourier domain 
\be \label{equ_sim} 
\tilde{q}_j  = \tilde{g}_j \tilde{h}_j, \quad  \nu_j = 1/\tilde{h}_j, \   j=1, \cdots, n,
\ee 
where, for any function $t$, we denote its $j$-th Fourier coefficient  by $\tilde{t}_j$.
The periodic Fourier basis   serves as the eigenbasis for this operator.

We carried out simulations with the periodized versions of the following two kernels  
\be\label{kernels}
g_1 (x)= 0.5\, \exp(-\lam |x|), \quad g_2 (x)=\exp( -\lam x^2/2)
\ee
where $g_1(x)$ corresponds to the case of $\alpha = \beta = 0, \gamma =2$ while 
$g_2(x)$ corresponds to $\alpha \propto 1/\lam$, $\beta = 2$  in \eqref{nu_j_cond}.
Hence, the problem is moderately ill-posed with $g_1$ and  severely ill-posed with $g_2$.
In addition, recovery of the solution becomes easier as $\lam$ grows.

Although we carried out simulations for a much wider sets of parameters, here we report the results for 
two  series of simulations with $n = 256$, $M=60$ and  $K=4$. 
In the first batch, we considered a set of smooth spatially homogeneous test functions 
\be    \label{smooth}
 l_1(x) = \sin(4\pi x),\ 
l_2(x)= \sin(4\pi(x-1/16)),\ 
l_3 (x) = \lkr x-0.5 \rkr^2,\ 
l_4 (x) = \lkr x-0.5\rkr^4,
\ee
coefficients of which follow the assumption \eqref{coef_cond}. For this set, we used 
Fourier basis $\phi_j$, $j=1, \cdots, n$, that diagonalizes the problem. Moreover, since 
the functions are spatially homogeneous, they can be well estimated when the same set $J$ 
of nonzero coefficients is used for all four functions.
In the second round, we   expanded our study to the set of spatially inhomogeneous functions
\be   \label{nonsmooth}
l_1(x) = l_B (x),\ 
l_2 (x) = l_W (x),\ 
l_3 (x) = l_P (x),\ 
l_4 (x) = |x-0.5|
\ee
where $l_B (x)$, $l_W (x)$ and $l_P(x)$ are the {\it blip, wave} and {\it parabolas} 
introduced by Donoho and  Johnstone \cite{donjohn}. In this case, Fourier basis does not allow accurate estimation,
hence, we used the Daubechies 8 wavelet basis as $\phi_j$, $j=1, \cdots, n$,
for which conditions  \eqref{Don1} and \eqref{Don3}  hold with $\nu_j$ given in \eqref{equ_sim} 
(see, e.g., \cite{abram}).  Although the second example does not follow our assumptions, 
it shows that our conclusions are true even in the situation when those assumptions are violated. 
In particular,  we used a different set of nonzero coefficients $J_k$
for $l_k$, $k=1,\ldots,4$,  for the functions in \eqref{nonsmooth}.
%
%
%
%
\begin{table}
\begin{center}
\begin{tabular}{|c|c |c | c| c | c |}
\hline \hline  
\multicolumn{6}{|c|}{$\lam=7$}\\
\hline 
  & \multicolumn{2} {c|} { Clustering Before} & \multicolumn{2} {c|}{ Clustering After} &{No Clustering }  \\
    \cline{2-5}
 {}   &{ Error} & { Miss-rate} &{ Error}  & { Miss-rate} &{ }  \\
    \hline \hline
$SNR=3$     &0.0365(0.0262) &0.0090 &0.0554(0.0083) &0.0068  &0.0556(0.0001) \\
\hline
$SNR=5$       &0.0270(0.0015) &0.0000 &0.0419(0.0095) &0.0070  &0.0423(0.0001)\\
\hline
$SNR=7$       &0.0250(0.0084) &0.0031 &0.0405(0.0000) &0.0000  &0.0414(0.0000)\\
\hline \hline
\multicolumn{6}{|c|}{$\lam=5$}\\
\hline 
$SNR=3$       &0.0377(0.0038) &0.0000 &0.0549(0.0059) &0.0033  &0.0567(0.0002)\\
\hline
$SNR=5$      &0.0317(0.0016) &0.0000 &0.0542(0.0000) &0.0000  &0.0551(0.0000)\\
\hline
$SNR=7$      &0.0269(0.0014) &0.0000 &0.0406(0.0000) &0.0000  &0.0421(0.0001)\\
\hline \hline
\multicolumn{6}{|c|}{$\lam=3$}\\
\hline
$SNR=3$      &0.0498(0.0398) &0.0106 &0.0810(0.0217) &0.0133  &0.0788(0.0001)\\
\hline
$SNR=5$      &0.0409(0.0237) &0.0036 &0.0543(0.0000) &0.0000  &0.0565(0.0002)\\
\hline
$SNR=7$      &0.0350(0.0241) &0.0026 &0.0542(0.0000) &0.0000  &0.0554(0.0001)\\
\hline
\end{tabular}
\end{center}
\caption{Estimation and clustering errors for the ``Clustering before'', ``Clustering after'' and  
``No clustering'' scenarios averaged over 100 simulation runs (the standard deviations of the means are in parentheses). 
Results for the set of  functions   \eqref{smooth} with the $g_1 (x)$ kernel in \eqref{kernels} and the same set of 
nonzero coefficients for all functions. }
\label{Table1}
\end{table}
%
%
%
%
%
%
We sampled the test functions on the equispaced grid on the interval $[0,1]$ and  scaled them to have norms $\sqrt{n}$, obtaining 
$h_k=c_k l_k$ where $c_k=\sqrt{n}/\|l_k\|$, $k=1,\ldots, 4$.     
Note that, while  the  functions in Set 1 \eqref{smooth} are simpler and easier to recover, 
they are less distinct and harder to cluster since $l_1$ is   similar to  $l_2$ and $l_3$ is similar to $l_4$.
On the other hand, while it is easier to distinguish between images of functions  in  Set 2 \eqref{nonsmooth},
they  are more difficult to estimate.   
For each of the test functions $h_k$, $k=1, \cdots, K$, we evaluated $u_k=(Ah)_k$, and sampled those functions on the 
 grid   of $n$ equispaced points $j/n$, $j=1, \cdots, n$, 
on the interval $[0,1]$, obtaining vectors $\bh_k$ and $\bu_k, k=1, \cdots, K$.
Furthermore, we generated a  clustering function $z: M \to K$ that places $M$ objects into $K$ classes,
$M/K$ into each class at random. We obtained the true matrices $\bF, \bQ \in \RR^{n \times M}$ with the 
columns   $\bh_{z(m)}$ and  $\bu_{z(m)}$, $m=1, \cdots, M$, respectively.
Finally, we generated data  $\bX$ by adding independent Gaussian noise with the 
standard deviation $\sig$ to every element in $\bQ$.
We found $\sig$ by fixing the Signal-to-Noise Ratio (SNR) and choosing 
$\sig = \std(\bF)/SNR$, where $\std(\bF)$ is the standard
deviation of the matrix   $\bF$ reshaped as a vector.
In what follows, we  considered  several noise scenarios: SNR = 3,  5   
and   7  for $g_1$ and SNR = 5, 7, and 10 for  $g_2$.  
In our study we treat $K$ as known and compare the estimators  where clustering was carried out at pre-processing level 
(``Clustering before'') to the estimators where clustering was done post-estimation (``Clustering after'') 
and estimators  without clustering (``No clustering''). 
 

\begin{table}
\begin{center}
\begin{tabular}{|c|c |c | c| c | c |}
\hline \hline
\multicolumn{6}{|c|}{$\lam=15$}\\
\hline 
  & \multicolumn{2} {c|} { Clustering Before} & \multicolumn{2} {c|}{ Clustering After} &{No Clustering }  \\
    \cline{2-5}
 {}   &{ Error} & { Miss-rate} &{ Error}  & { Miss-rate} &{ }  \\
    \hline 
$SNR=5$      &0.1568(0.0684) &0.0623 &0.1258(0.0175) &0.0071  & 0.1252(0.0002) \\

\hline
$SNR=7$       &0.1516(0.0640) &0.0521 &0.1252(0.0063) &0.0180  &0.1245(0.0001)\\

\hline
$SNR=10$       &0.1307(0.0342) &0.0128 &0.1237(0.0000) &0.0000  &0.1241(0.0000)\\

\hline \hline
\multicolumn{6}{|c|}{$\lam=12$}\\

\hline
$SNR=5$      &0.2336(0.0770) &0.1601 &0.1609(0.0173) &0.0398  &0.1659(0.0034)\\

\hline
$SNR=7$     &0.2186(0.0759) &0.1303 &0.1602(0.0080) &0.0413  &0.1620(0.0025)\\

\hline
$SNR=10$      &0.1938(0.0660) &0.0758 &0.1583(0.0058) &0.0211  &0.1592(0.0017)\\

\hline \hline
\multicolumn{6}{|c|}{$\lam=10$}\\
\hline

$SNR=5$      &0.5419(0.0707) &0.2513 &0.7933(0.1005) &0.2796  &0.7448(0.0000)\\

\hline
$SNR=7$      &0.5196(0.0128) &0.2331 &0.7678(0.0733) &0.2693  &0.7448(0.0000)\\

\hline
$SNR=10$      &0.5078(0.0212) &0.1715 &0.4853(0.0084) &0.0430  &0.4849(0.0037)\\
\hline
\end{tabular}
\end{center}

\caption{Estimation and clustering errors for the ``Clustering before'', ``Clustering after'' and  
``No clustering'' scenarios averaged over 100 simulation runs (the standard deviations of the means are in parentheses). 
Results for the set of  functions   \eqref{smooth} with the $g_2 (x)$ kernel in \eqref{kernels} and the same set of 
nonzero coefficients for all functions. }
\label{Table2}
\end{table}


For the ``Clustering before'' setting, we applied  clustering directly to the elements of matrix $\bY$. 
As it follows from equation \eqref{opt_prob2}, the matrix $\hat{\bZ} \in \calM (M,K)$ which minimizes the objective function
is a solution of the $K$-means clustering problem. Subsequently, we found matrix $\bPi_{\hat{\bZ},K}$
and, following equation \eqref{eq:opt_bG}, estimated $\bG_*$ by $\widehat{\bG} = \bW_{ \hat{J}} \bUp \bY \bPi_{\hat{\bZ}, K}$.
For the  set of functions \eqref{smooth}, the set $\hat{J}$ was obtained 
by applying hard thresholding to the rows of the matrix $\bUp \bY \bPi_{\hat{\bZ}, K}$,
while for the  set of functions \eqref{nonsmooth}, we applied hard 
hard thresholding to each of the elements of the matrix $\bUp \bY \bPi_{\hat{\bZ}, K}$.
Finally, the estimator $\widehat{\bF}$ of the matrix $\bF_*$ is obtained by 
applying the inverse Fourier  transform (in the case of the functions in \eqref{smooth}) 
or the inverse   wavelet transform (in the case of the functions in \eqref{nonsmooth})  
to the columns of  matrix $\widehat{\bG}$.
For the ``Clustering after'' setting, we first constructed the ``No clustering'' estimator $\check{\bG}$ of  matrix  $\bG_*$ 
by thresholding elements of the columns of the matrix $\bUp \bY$ in equation \eqref{main_eq1}, and then obtained 
the estimator $\check{\bF}$ of  matrix  $\bF_*$ by applying the inverse Fourier or wavelet transform to the columns of  matrix $\check{\bG}$.
Finally, the ``Clustering after''  estimator of $\bF_*$ is obtained by applying the $K$-means clustering procedure 
to the columns of matrix $\check{\bF}$.


\begin{table}
\begin{center}
\begin{tabular}{|c|c |c | c| c | c |}
\hline \hline
\multicolumn{6}{|c|}{$\lam=7$}\\
\hline
  & \multicolumn{2} {c|} { Clustering Before} & \multicolumn{2} {c|}{ Clustering After} &{No Clustering }  \\
    \cline{2-5}
 {}   &{ Error} & { Miss-rate} &{ Error}  & { Miss-rate} &{ }  \\
    \hline
$SNR=3$     &0.1364(0.0055) &0.0000 &0.2650(0.0609) &0.0250  &0.2810(0.0056) \\
\hline
$SNR=5$      &0.1190(0.0039) &0.0000  &0.2187(0.0661) &0.0180 &0.2470(0.0030)\\
\hline
$SNR=7$      &0.1033(0.0058) &0.0000  &0.1700(0.0599) &0.0205 &0.2004(0.0039)\\
\hline \hline
\multicolumn{6}{|c|}{$\lam=5$}\\
\hline
$SNR=3$      &0.1480(0.0077) &0.0000 &0.2845(0.1095) &0.0731 &0.3569(0.0091)\\
\hline
$SNR=5$     &0.1186(0.0053) &0.0000  &0.2322(0.1117) &0.0610  &0.2719(0.0040)\\
\hline
$SNR=7$     &0.1026(0.0045) &0.0000  &0.1632(0.0656) &0.0221  &0.2169(0.0042)\\
\hline \hline
%
\multicolumn{6}{|c|}{$\lam=3$}\\
\hline
$SNR=3$     &0.1806(0.0110) &0.0000  &0.2932(0.1309) &0.1010 &0.4831(0.0092)\\
\hline
$SNR=5$     &0.1442(0.0061) &0.0000  &0.2326(0.1199) &0.0690 &0.3250(0.0053)\\
\hline
$SNR=7$     &0.1310(0.0047) &0.0000 &0.2149(0.1207) &0.0718 &0.2542(0.0042)\\
\hline
\end{tabular}
\end{center}
\caption{Estimation and clustering errors for the ``Clustering before'', ``Clustering after'' and  
``No clustering'' scenarios averaged over 100 simulation runs (the standard deviations of the means are in parentheses). 
Results for the set of  functions   \eqref{nonsmooth}  with the $g_1 (x)$ kernel in \eqref{kernels} and unique set of 
nonzero coefficients for each  of the functions. }
\label{Table3}
\end{table}




\begin{table}
\begin{center}
\begin{tabular}{|c|c |c | c| c | c |}
\hline \hline
\multicolumn{6}{|c|}{$\lam=15$}\\
\hline 
  & \multicolumn{2} {c|} { Clustering Before} & \multicolumn{2} {c|}{ Clustering After} &{No Clustering }  \\
    \cline{2-5}
 {}   &{ Error} & { Miss-rate} &{ Error}  & { Miss-rate} &{ }  \\
    \hline
 
 
$SNR=5$     &0.3709(0.0000) &0.0000 &0.3709(0.0000) &0.0000  &0.3714(0.0001) \\

\hline
$SNR=7$      &0.3708(0.0000) &0.0000 &0.3708(0.0000) &0.0000  &0.3711(0.0000)\\

\hline
$SNR=10$      &0.3708(0.0000) &0.0000  &0.3708(0.0000) &0.0000 &0.3710(0.0000)\\

\hline \hline
\multicolumn{6}{|c|}{$\lam=12$}\\

\hline
$SNR=5$      &0.3768(0.0009) &0.0000 &0.3768(0.0009) &0.0000  &0.3810(0.0011)\\

\hline
$SNR=7$     &0.3766(0.0006) &0.0000  &0.3780(0.0137) &0.0036 &0.3787(0.0006)\\

\hline
$SNR=10$     &0.3765(0.0004) &0.0000  &0.3785(0.0202) &0.0035 &0.3776(0.0004)\\

\hline \hline
\multicolumn{6}{|c|}{$\lam=10$}\\
\hline
$SNR=5$     &0.4876(0.0049) &0.0000  &0.4933(0.0294) &0.0141  &0.4940(0.0050)\\

\hline
$SNR=7$     &0.4869(0.0035) &0.0000  &0.4869(0.0035) &0.0000 &0.4903(0.0035)\\

\hline
$SNR=10$     &0.4872(0.0027) &0.0000  &0.4872(0.0027) &0.0000 &0.4888(0.0027)\\
\hline
\end{tabular}
\end{center}
\caption{Estimation and clustering errors for the ``Clustering before'', ``Clustering after'' and  
``No clustering'' scenarios averaged over 100 simulation runs (the standard deviations of the means are in parentheses). 
Results for the set of  functions   \eqref{nonsmooth}  with the $g_2 (x)$ kernel in \eqref{kernels} and unique set of 
nonzero coefficients for each  of the functions. }
\label{Table4}
\end{table}



Tables 1--4 report simulations results for the three clustering scenarios above  
(``Clustering before'', ``Clustering after'' and ``No clustering''), 
for each of the sets of test functions in \eqref{smooth} and  \eqref{nonsmooth}  and 
for each of the two kernels in \eqref{kernels}  with various values of $\lam$. 
In the Tables, we display the accuracies of the three estimators where the precision of an   estimator $\widehat{\bF}$
is measured by  the  Frobenius norms of its error   
\be\label{relerror}
\Del = \Del(\widehat{\bF}) = \|\widehat{\bF} - \bF \|_F/\sqrt{Mn}.  
\ee
In addition,   we report the proportion of  erroneously clustered nodes (``Miss-rate'')
for the ``Clustering before''and the ``Clustering after''   estimators.

We ought to point out that the ``Clustering before'' estimation procedure is much more computationally 
efficient since it does not require to recover $M$ unknown functions separately which is necessary for the 
``Clustering after'' and ``No clustering'' procedures.



\begin{figure} 
\[\includegraphics[height=4.0cm]{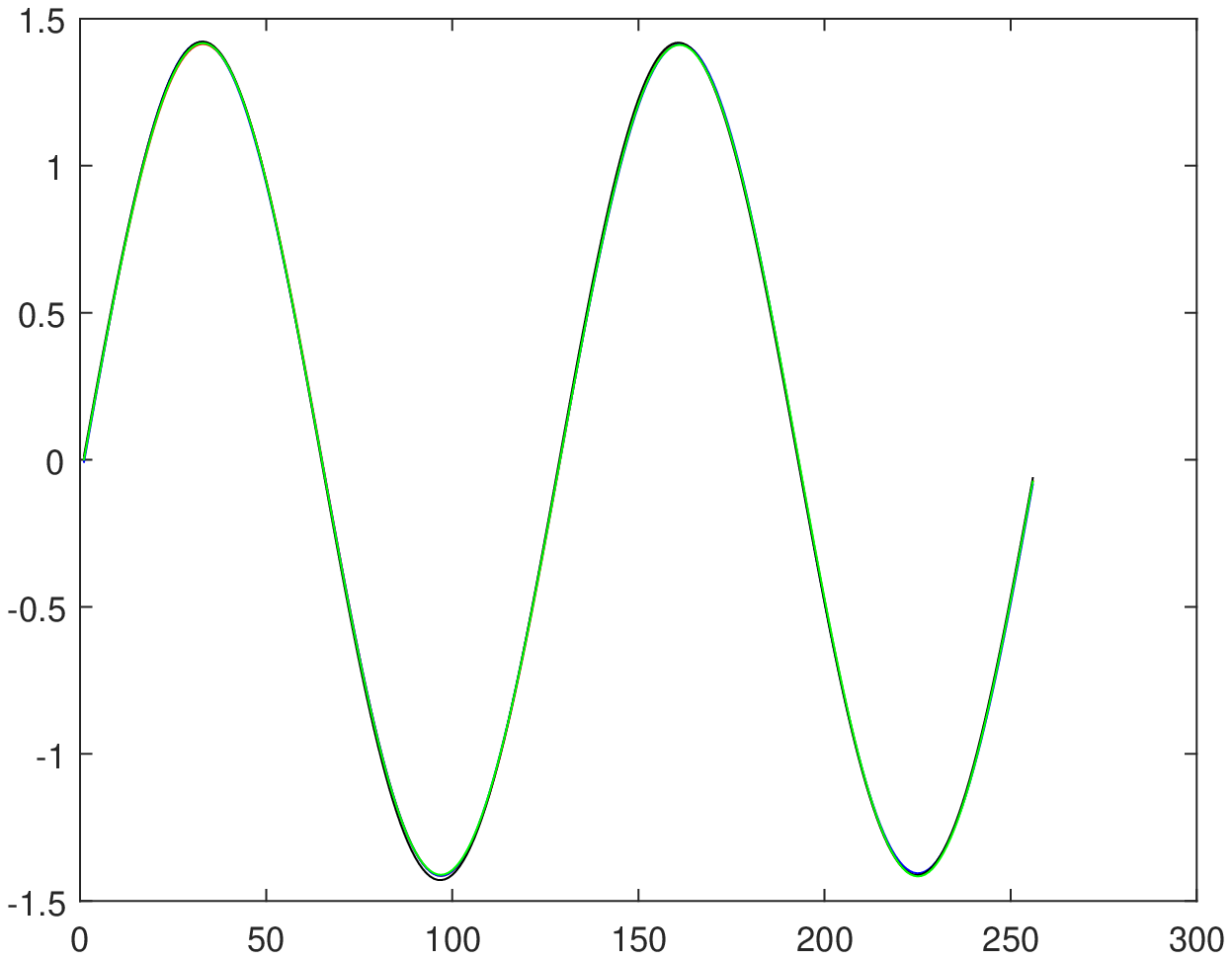} \hspace{2mm} 
\includegraphics[height=4.0cm]{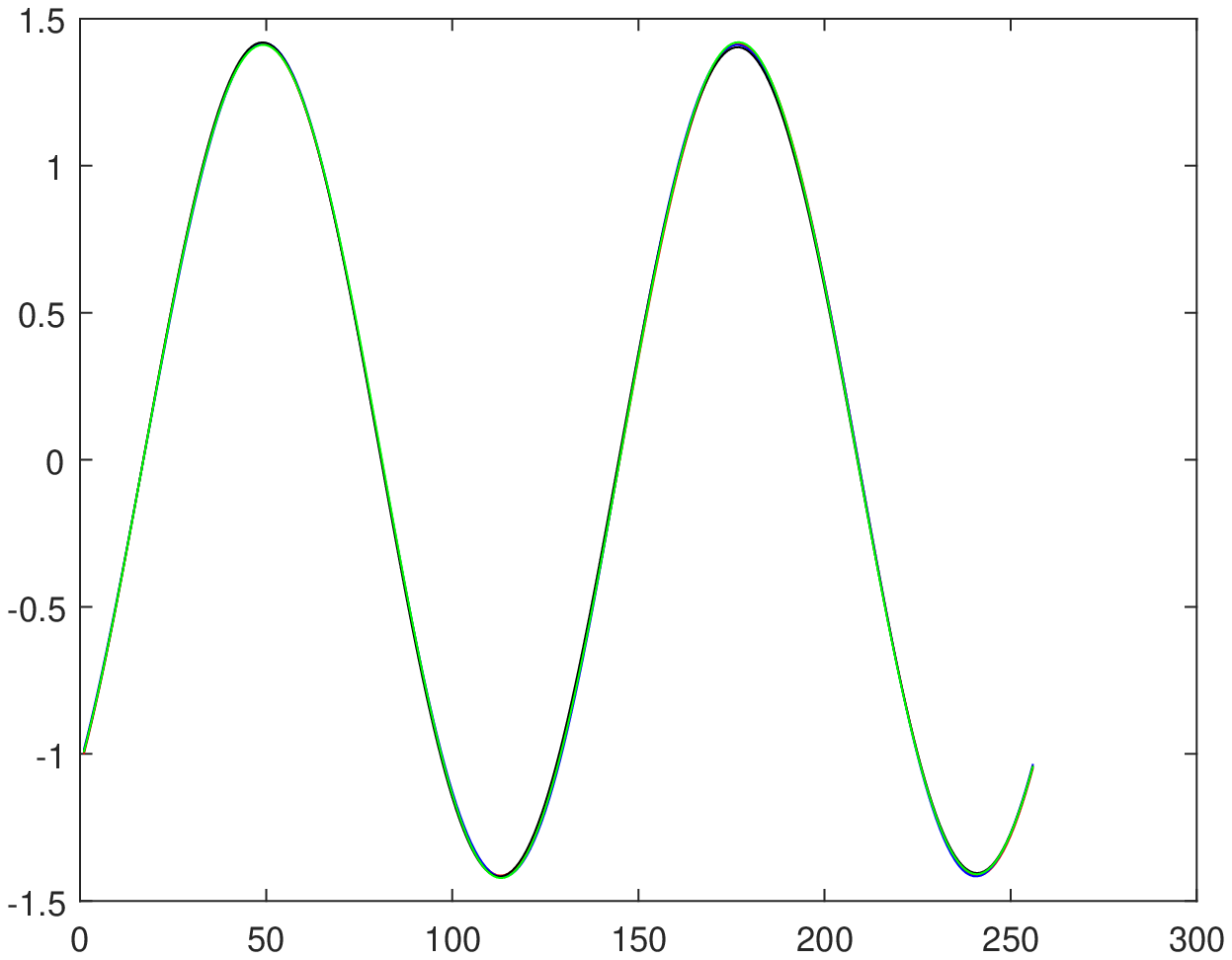} \] 
\[\includegraphics[height=4.0cm]{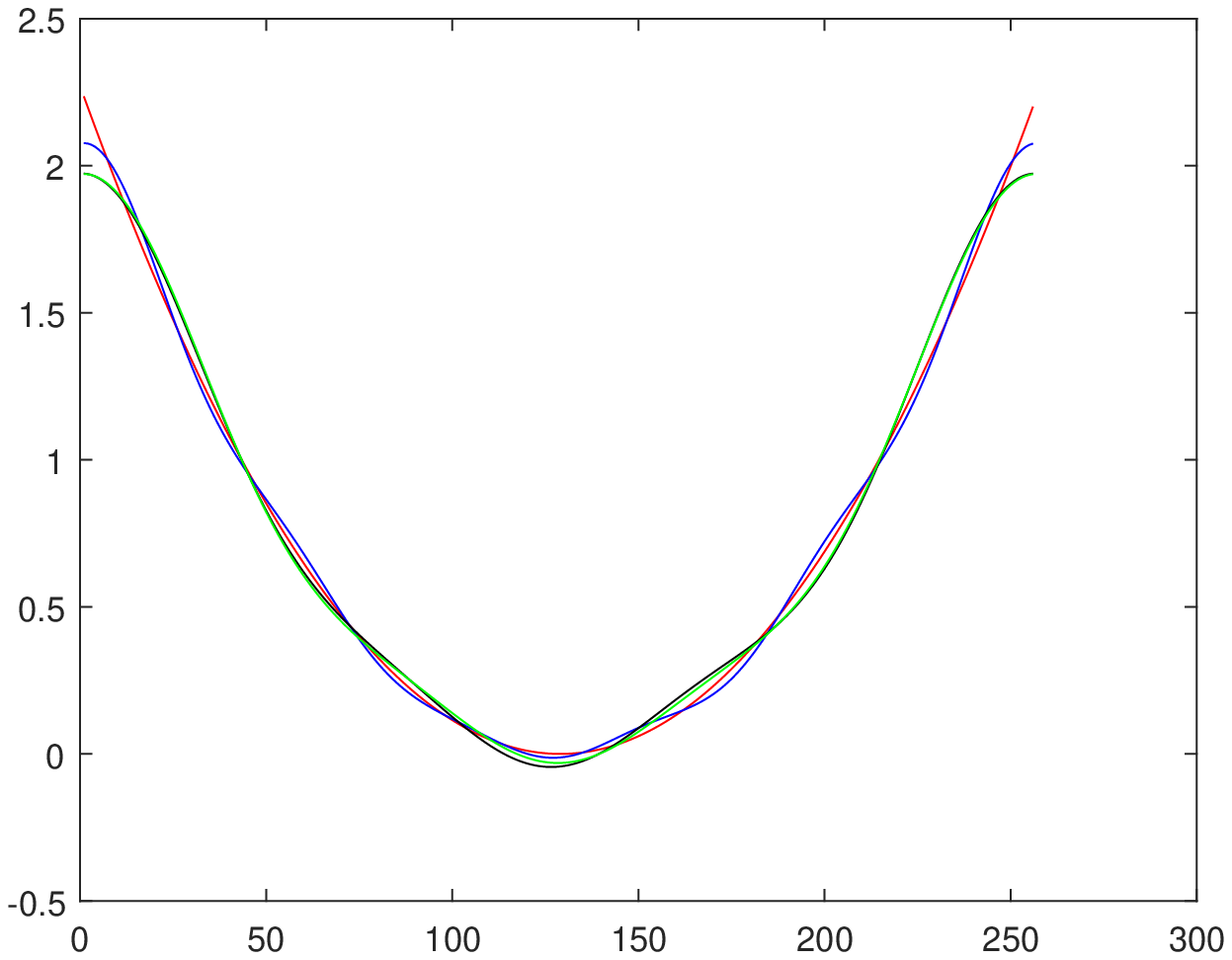} \hspace{2mm}
\includegraphics[height=4.0cm]{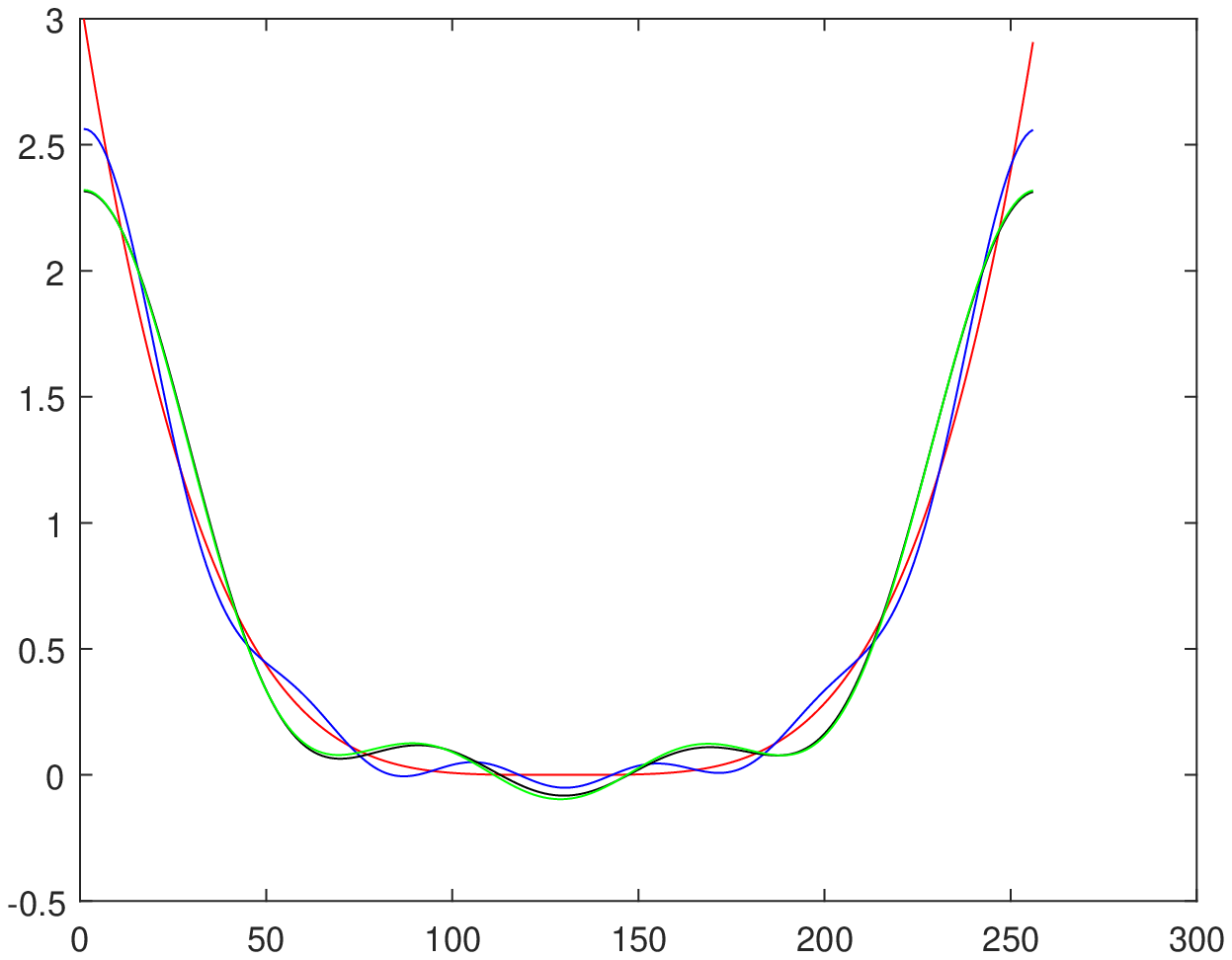}  
\]
\caption{True functions (red) and their estimators: ``Clustering before'' (blue),
``Clustering after'' (green) and ``No clustering'' (black). Results for the  
functions in \eqref{smooth} and the kernel $g_1$ in \eqref{kernels} with $\lam=3$ and SNR=3.
Top row:  $h_1$ (left), $h_2$ (right). Bottom row:   $h_3$ (left), $h_4$ (right).   
\label{fig1}}
\end{figure}


\begin{figure} 
\[\includegraphics[height=4.0cm]{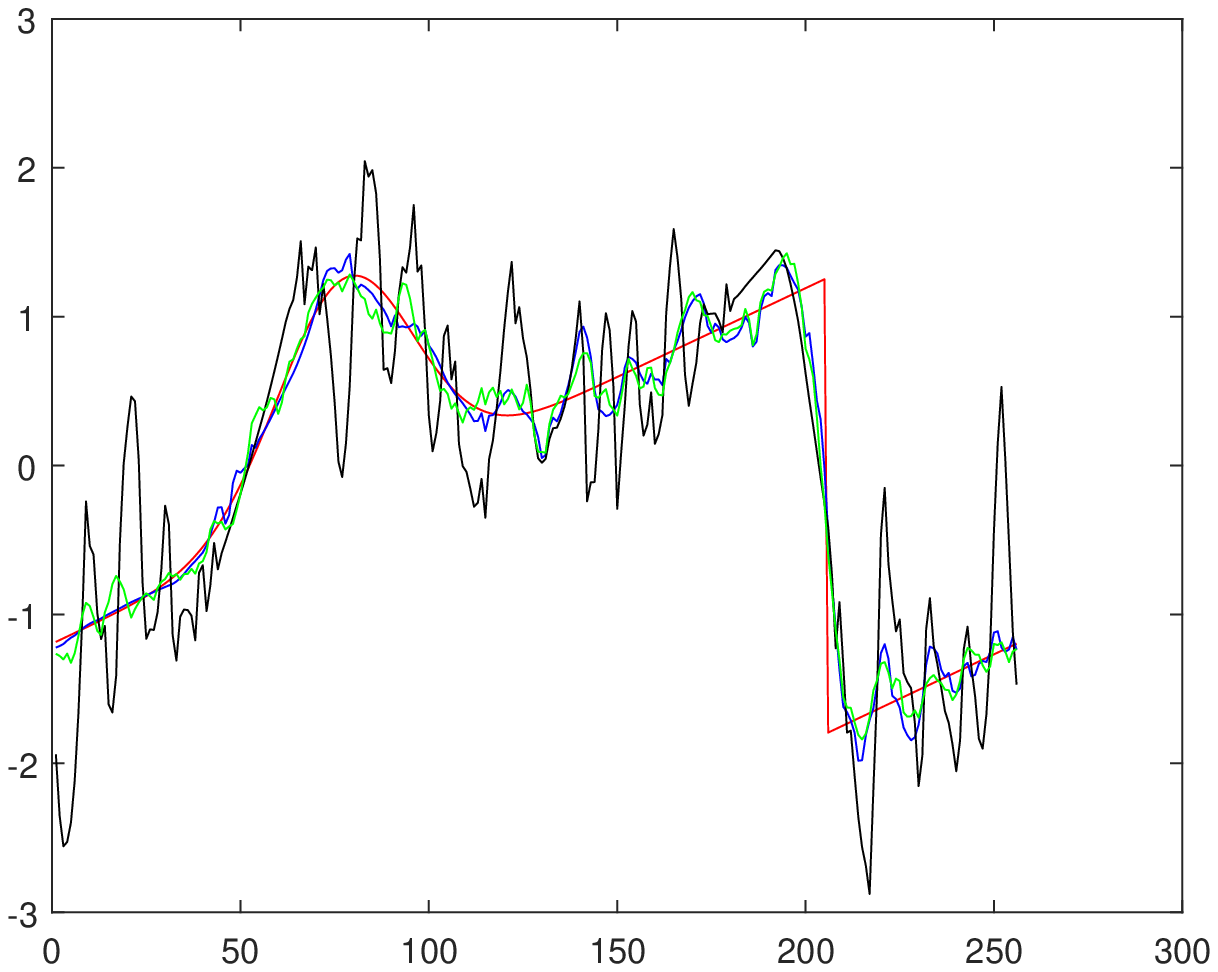} \hspace{2mm} 
\includegraphics[height=4.0cm]{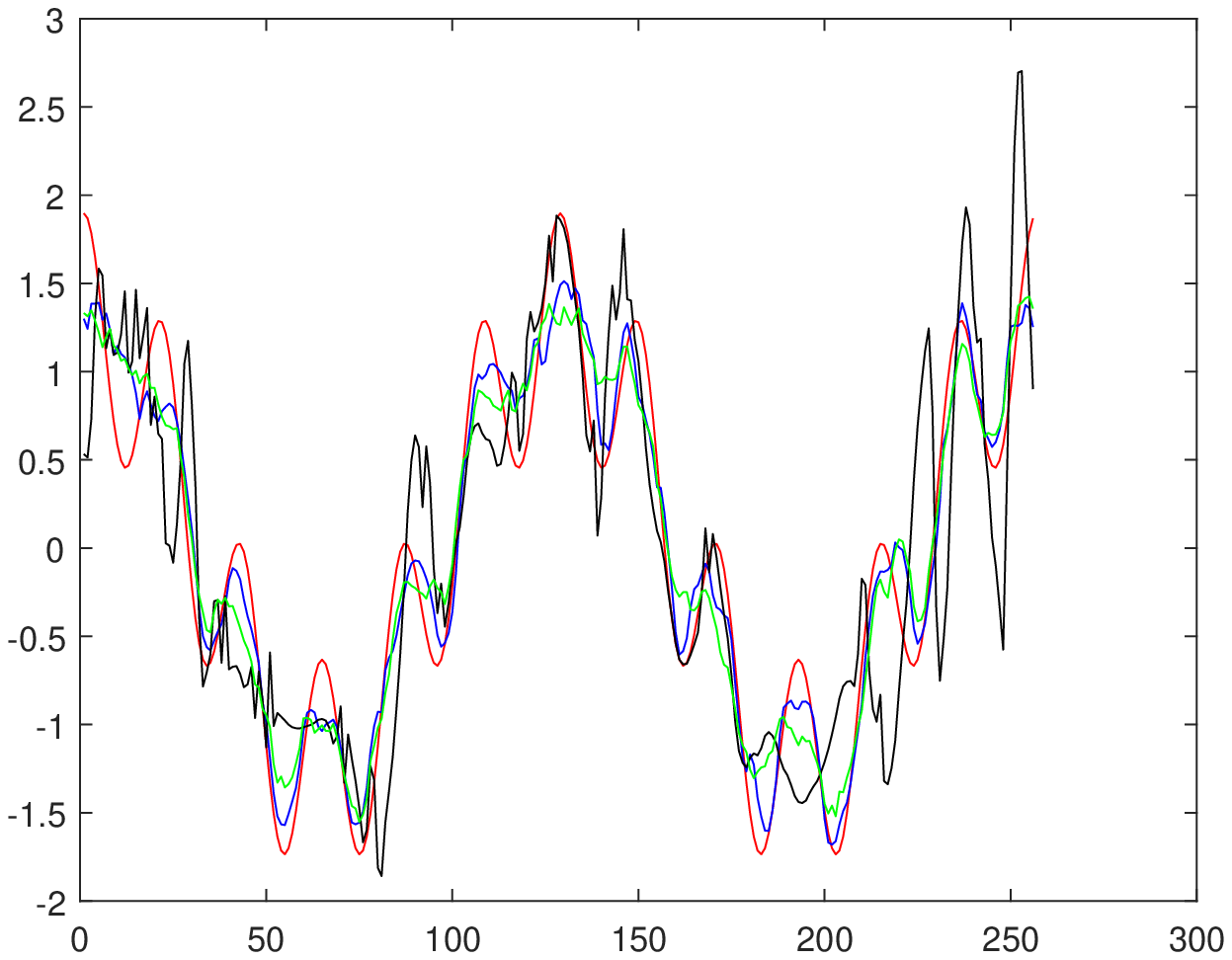} \] 
\[\includegraphics[height=4.0cm]{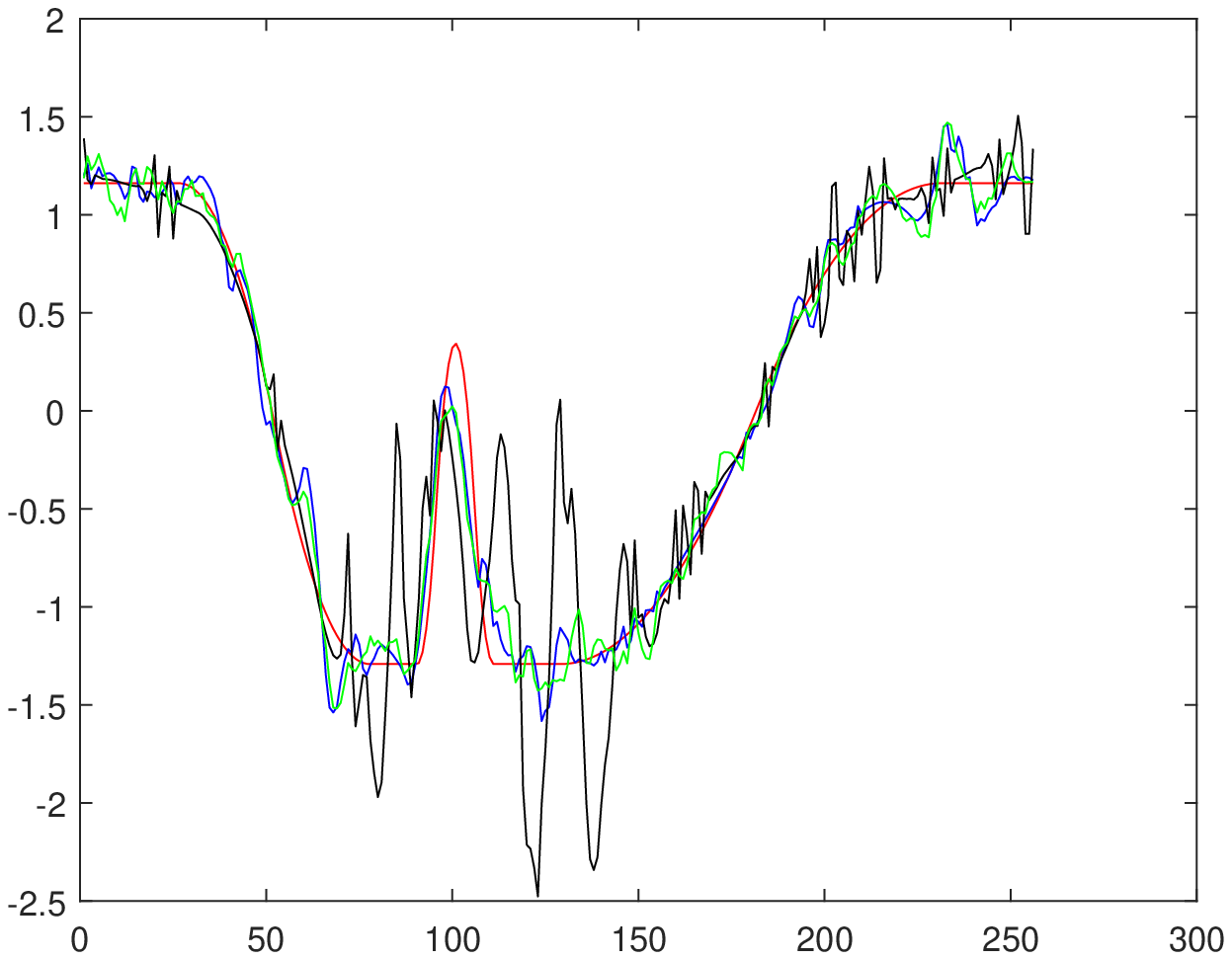} \hspace{2mm}
\includegraphics[height=4.0cm]{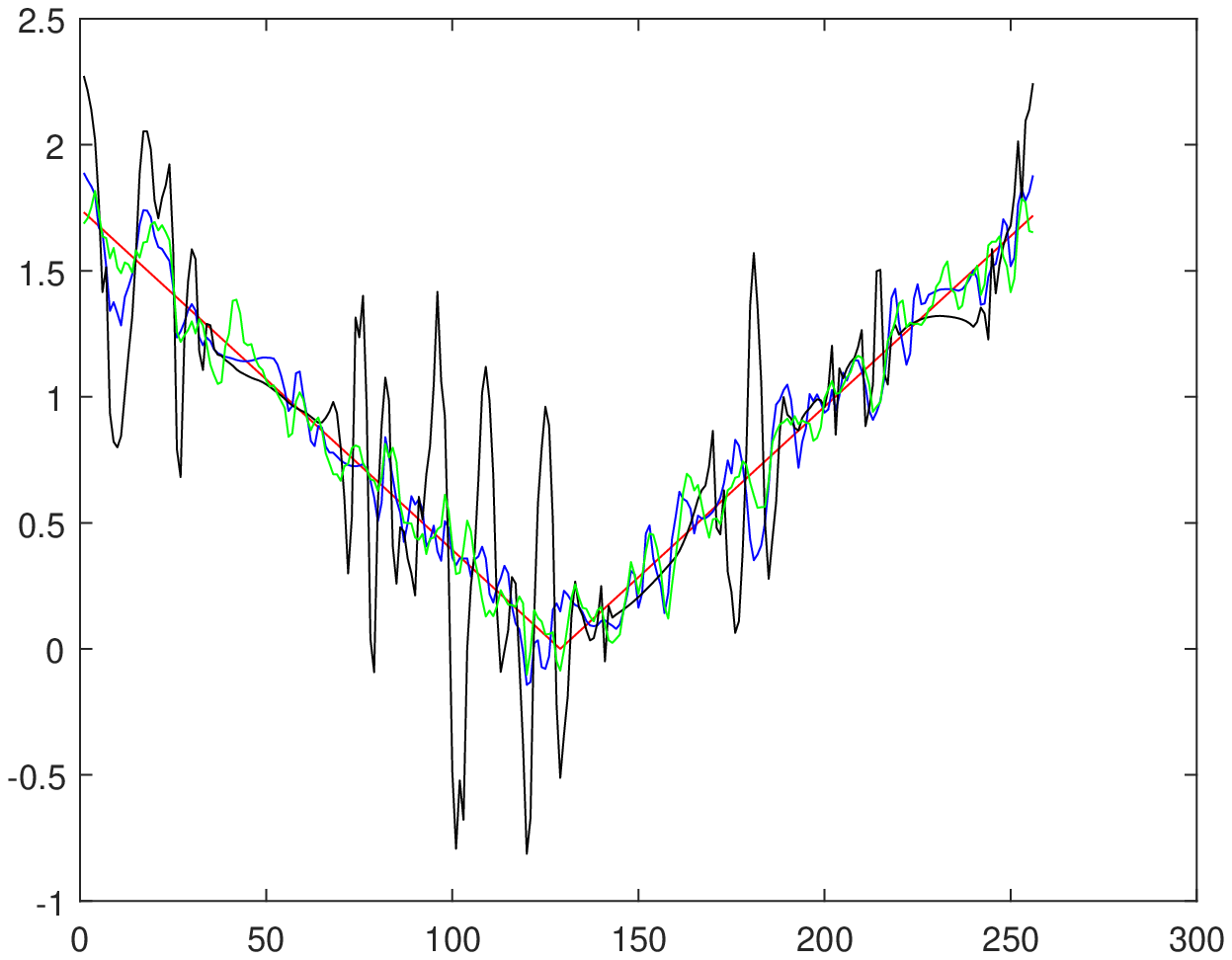}  
\]
\caption{True functions (red) and their estimators: ``Clustering before'' (blue),
``Clustering after'' (green) and ``No clustering'' (black). Results for the  
functions in \eqref{nonsmooth} and the kernel $g_1$ in \eqref{kernels} with $\lam=3$ and SNR=3.
Top row:  $h_1$ (left), $h_2$ (right). Bottom row:   $h_3$ (left), $h_4$ (right).  
\label{fig2}}
\end{figure}


 
\section {Conclusions }
\label{sec:discussion}
\setcounter{equation}{0}

In this paper, we investigate theoretically and via a limited simulation study, the effect of clustering
on the accuracy of recovery in   ill-posed linear inverse problems. 
As we have stated earlier, in many applications leading to such  problems, clustering is carried 
out at a pre-processing step and later is totally forgotten when it comes to error evaluation.
Our main objective has been to evaluate what effect  clustering at the pre-processing step   has 
 on the precision of the resulting estimators.  

It appears that benefits of pre-clustering depend significantly on the nature of the inverse problem 
at hand. If the problem is moderately ill-posed (kernel $g_1$ in \eqref{kernels}, $\alpha= \beta =0$), 
then, as    Corollary \ref{cor:comparison}  shows, the ``Clustering Before'' estimator has 
asymptotically smaller errors than the ``No Clustering'' estimator when the number of functions and the sample size  grow.
Tables~1 and 2, corresponding to this case, confirm that, for the finite number of functions and moderate sample size,  the 
``Clustering before'' procedure delivers better precision than the ``Clustering after''   and ``No clustering'' 
techniques. Furthermore, the ``Clustering before''  estimation has profound computational benefits since one needs to 
recover $K$ unknown functions instead of $M$.   
Moreover, the advantages of clustering at pre-processing step become more prominent when the problem is less ill-posed
(larger  $\lam$). Indeed, in the case when the problem is not ill-posed ($\alpha = \beta = \gamma =0$  in \eqref{nu_j_cond}),
as findings of Klopp \etal \cite{klopp} show, clustering always improves estimation precision.

The situation changes drastically if the inverse problem is severely ill-posed 
(kernel $g_2$ in \eqref{kernels}, $\alpha >0,  \beta >0$). Our theoretical results indicate that clustering, in this case, does not improve 
the estimation precision as the number of functions and the sample size  grow. These findings are consistent with the simulation study.
In the case of functions in \eqref{nonsmooth}, Table~4 implies that the precisions of all three methodology are approximately 
the same, and the estimation errors are high even when clustering errors are small or zero.
This is due to the fact  that  the reduction in the noise level due to clustering is not sufficient to counteract 
the ill-posedness of the problem and, thus,  it does not lead to a meaningful improvement in estimation accuracy.
Table~3, that reports on the simulations with   functions in \eqref{smooth}, 
presents an even more grim picture. Since functions in the set \eqref{smooth} resemble  each other to start with and 
convolutions with the kernel $g_2$ make them to appear  even more similar, ``Clustering before'' procedure  leads 
to relatively high clustering errors that, in turn, produce higher estimation errors than the ``Clustering after''   and ``No clustering''
techniques. 
\\
 
In conclusion, clustering at the pre-processing step is beneficial when the problem is moderately ill-posed.
It should be applied with extreme care when the problem is severely ill-posed.


\section*{Acknowledgments}

 Marianna Pensky  and  Rasika Rajapakshage were  partially supported by National Science Foundation
 (NSF), grants     DMS-1407475 and DMS-1712977.


\input{RP_Proofs_March22_2020}


\input{References}

\end{document}

%% file: RP_Proofs_March22_2020.tex


\section{Proofs}
\label{sec:proofs}
\setcounter{equation}{0}


\subsection{Proof  of the oracle inequality}

\noindent
{\bf Proof of Theorem~\ref{th:upper_bound}. } \ 
The proof of the  inequality \eqref{eq:upper_bound} is based on the standard techniques for   proofs of oracle inequalities.
We use optimization problem \eqref{opt_prob1} to present the left-hand side as a sum of the error of any estimator plus the random error term 
followed by the difference between the penalty terms. Later on, we upper-bound the random error term for any  number of classes $K$, 
any clustering matrix $Z$ and any set $J$. After that, we take a union bound over all possible $K$, $Z$ and $J$ 
to obtain an upper bound for the probability that the error exceeds certain threshold. 
The novelty of the proof lies in the fact that we are using vectorization of the model which allows us to 
attain the upper bounds.

Note that it follows from the optimization problem  \eqref{opt_prob1}  that for any fixed 
 $\bG,\bZ,J$ and $K$ one has   
\bes
\|\widehat{\bG}\|_F ^2-2\Tr (\bY^T\bUp \widehat{\bG} \bPi_{\hat{\bZ},\hat{K}}) +\Pen(\hat{J},\hat{K})
\leq  \|\bG\|_F ^2-2\Tr (\bY^T\bUp \bG\bPi_{\bZ,K}) +\Pen(J,K).
\ees
Then, adding and subtracting $\bG_*$, we obtain
\begin{align*}
& \|\widehat{\bG}-\bG_* \|_F^2 + \|\bG_* \|_F^2 + 2 \Tr((\widehat{\bG}-\bG_*)^T \bG_*)
- 2 \Tr(\bY^T\bUp \widehat{\bG} \bPi_{\hat{\bZ},\hat{K}}) +\Pen(\hat{J},\hat{K})
\leq  \\
& \|\bG-\bG_* \|_F^2 + \|\bG_* \|_F^2 + 2 \Tr((\bG -\bG_*)^T \bG_*) 
-2\Tr (\bY^T\bUp \bG\bPi_{\bZ,K}) +\Pen(J,K).
\end{align*}
Combine the trace product terms and recall  that, due to equation \eqref{main_seq},
$\bY = \bUp^{-1} \bG_* + \del \bE$. Hence, the last inequality yields
\be\label{Lasso_1}
\|\widehat{\bG}-\bG_*\|_F^2  \leq  \|\bG-\bG_*\|_F^2+  {2\,\del}\, \Tr [\bE^T \bUp ( \widehat{\bG}  -\bG )] 
+\Pen(J,K) -\Pen(\hat{J},\hat{K})
\ee
We choose $\bG=\bW_J \bG_* \bPi_{\bZ,K}$ and, 
in order to analyze the cross term $\Tr [ \bE^T \bUp ( \widehat{\bG}  -\bG )]$, we use vectorization of the model.
For this purpose, we choose $\bS$ such that $\bSig = \bS \bS^T$    and denote
\begin{align} 
& \bPi_{\hat{\bZ},\hat{K},\hat{J}}=(\bPi_{\hat{\bZ},\hat{\bK}} \otimes \bW_{\hJ}),\quad  \bPi_{\bZ,K,J}=(\bPi_{\bZ,K} \otimes \bW_{J})
\label{vec_notations1}\\
& \hat{\bg}=\vect  (\widehat{\bG}), \quad \bg=\vect(\bG),\quad \beps = \vect(\bE), \quad  \bGa=(\bI_M \otimes \bUp), \quad  
\boeta = (\bI_M \otimes \bS^{-1}) \beps.  
\label{vec_notations2} 
\end{align}
By definition of the matrix-variate normal distribution  (Theorem 2.3.1 of 
Gupta and Nagar (2000)) and \eqref{eq:mom_matr_normal}, we derive that 
\be \label{error_vect}
\beps  \sim  N(0,\bSig  \otimes \bI_M) 
\ee
Then, $\EE(\boeta\boeta^T) =\bI_{n M}$, so that,  $\boeta \sim N(0,\bI_{n M})$, 
where $\beps$ is defined in \eqref{vec_notations2} and $\|\bS\|_{op} \leq C_\psi$.
Then, equation \eqref{main_eq1} can be re-written as 
\be\label{vecfun} 
\bGa \by =\bg_* +  \del \, \bGa \, (\bI_M \otimes \bS) \boeta. 
\ee
Observe  that by  Theorem 1.2.22 of Gupta and Nagar (2000), one has
\bes
\hat{\bg}=\vect (\bW_{\hat{J}}\bUp \bY \bPi_{\hat{\bZ},\hat{K}})= 
\bPi_{\hat{\bZ},\hat{K},\hat{J}}\bGa \by,\quad \bg=\bPi_{\bZ,K,J}\bg_*
\ees
and  $\Tr[\bE^T \bUp ( \widehat{\bG}  -\bG )] =\boeta^T(\bI_M \otimes \bS^T\bUp )
(\bPi_{\hat{\bZ},\hat{K},\hat{J}}\bGa \by-\bPi_{\bZ,K,J}\bg_*)$.
Now  \eqref{Lasso_1}  can be rewritten in a vector form as
\be\label{VLasso}
\|\hat{\bg}-\bg_*\| ^2  \leq  \|\bg-\bg_*\| ^2 + \Del
+\Pen(J,K) -\Pen(\hat{J},\hat{K})
\ee
where 
\be \label{eq:Del}
\Del= {2\,\del} \boeta^T(\bI_M \otimes \bS^T\bUp )
(\bPi_{\hat{\bZ},\hat{K},\hat{J}}\bGa \by-\bPi_{\bZ,K,J}\bg_*) = \Del_1+\Del_2
\ee 
with 
\be \label{Del1_2}
\Del_1= {2\,\del} \boeta^T(\bI_M \otimes \bS^T\bUp )(\bPi_{\hat{\bZ},\hat{K},\hat{J}}(\bGa \by-\bg_*)), \quad
\Del_2= {2\, \del} \boeta^T(\bI_M \otimes \bS^T\bUp )(\bPi_{\hat{\bZ},\hat{K},\hat{J}}-\bPi_{\bZ,K,J})\bg_*.
\ee 
 
Derivation of upper bounds for $\Del_1$ and $\Del_2$ is based on the following lemma.


\begin{lemma}\label{lem:lemma1}
Let $K,J$  be fixed,  $\hat{J}$ be an arbitrary random subset of $\{1, \ldots, n\}$
and $\hat{K}$ be a random integer between 1 and $M$. 
Let $\bZ \in \calM (M,K)$ and $\hbZ \in \calM (M,\hat{K})$ be a fixed and a random clustering matrix, respectively.
Denote the projection matrices on the column spaces of matrices $\bZ$ and $\hbZ$ by, respectively,   $\bPi_{\bZ,K}$ and $\bPi_{\hbZ,\hK}$. 
Let $\bS$ be a matrix with $\|\bS\|_{op} \leq C_\psi$ and $\boeta \sim N(0,\bI_{nM})$. 
Then, for any $\tau >0$, there exist  sets $\Om_{1\tau}$ and $\Om_{2\tau}$  with 
$\PP(\Om_{1\tau}) \geq 1 - \del^{\tau}$ and  $\PP(\Om_{2\tau}) \geq 1 - \del^{\tau}$ 
such that     
\be \label{eq1:lemma1}
 \|(\bPi_{\bZ,K}  \otimes (\bW_{J} \Up \bS)) \boeta \|^2   \leq 2K C_\psi ^2(\sum_{j \in J}\nu_j ^2) 
  + 3C_\psi ^2(\max_{j \in J} \nu_j ^2) \tau\, \ln (\del^{-1}) , \quad  \forall \om \in \Om_{1\tau};
\ee
\be \label{eq2:lemma1}
\begin{aligned}
 \|( \bPi_{\hbZ,\hK} & \otimes (\bW_{\hJ} \Up \bS)) \boeta  \|^2   \leq 2\hK C_\psi ^2(\sum_{j \in \hat{J}}\nu_j ^2)
 \\
& + 3C_\psi ^2(\max_{j \in \hat{J}} \nu_j ^2) \lfi M \ln \hK + |\hat{J}|\ln (ne/|\hat{J}|) +\ln (M n) +  \tau\, \ln (\del^{-1})   \rfi 
\quad  \forall \om \in \Om_{2\tau}.
\end{aligned}
\ee
Moreover, if  $J = \lfi 1, ..., L \rfi$ is fixed and $\hat{J} = \lfi 1, ..., \hat{L} \rfi$ 
for some random integer $\hat{L} \geq 1$, then%
\be \label{eq1:lemma1_new}
 \|(\bPi_{\bZ,K}  \otimes (\bW_{J} \Up \bS)) \boeta \|^2   \leq 2K C_\psi^2 \, \sum_{j=1}^L \nu_j^2  
  + 3C_\psi^2\,  \tau\, \ln (\del^{-1}) \, \nu_L^2, \quad  \forall \om \in \Om_{1\tau};
\ee
\be \label{eq2:lemma1_new}
\begin{aligned}
 \|( \bPi_{\hbZ,\hK} & \otimes (\bW_{\hJ} \Up \bS)) \boeta  \|^2   \leq 2 \hK C_\psi^2 \, \sum_{j=1}^L \nu_j^2  
 \\
& + 3C_\psi ^2\, \nu_L^2\,  \lfi M \ln \hK +    \ln (M n) +  \tau\, \ln (\del^{-1})   \rfi 
\quad  \forall \om \in \Om_{2\tau}.
\end{aligned}
\ee
\end{lemma}

\medskip


\noindent 
In what follows, we carry out only the proof of the upper bound \eqref{eq:upper_bound} that takes place for a generic
set $J$. The proof of the upper bound \eqref{eq:upper_bound_new} can be obtained from the proof below with minimal 
modifications.

Note that $\Del_1$ can be re-written as 
$\Del_1=  2\,\del^2\, \boeta^T(\bI_M \otimes \bS^T\bUp )(\bPi_{\hat{\bZ},\hat{K}} \otimes \bW_{\hJ})(\bI_M \otimes \bUp\bS )\boeta$.
Due to  $\bGa \by-\bg_* = \del \, \bGa \beps$  and  \eqref{vec_notations2}, we obtain
$\Del_1= 2\,\del^2\,   \|( \bPi_{\hat{\bZ},\hat{K}}  \otimes (\bW_{\hat{J}} \Up \bS)) \boeta \|^2$.
Therefore, by \fr{eq2:lemma1}, we obtain that for $\om \in \Om_{2\tau}$
\be \label{eq:Del1}
|\Del_1| \leq   {2\, \del^2\, C_\psi^2}\,  \lkv  2 \hK \, \sum_{j \in \hat{J}}\nu_j ^2 
 + 3 (\max_{j \in \hat{J}} \nu_j ^2) \lfi M \ln \hK + |\hat{J}|\ln (ne/|\hat{J}|)  +\ln (M n \del^{-\tau})  \rfi \rkv
\ee


In order to construct an upper bound for $\Del_2$, consider the following sets 
\be \label{setsJ}
\tJ = J \cup \hat{J},\quad J_1=J \cap \hat{J},\quad J_2= J^c \cap \hat{J}, \quad J_3= \hat{J}^c \cap J. 
\ee
The sets $J_1$, $J_2$ and $J_3$ are non-overlapping and $\tJ = J_1 \cup J_2 \cup J_3$.
Furthermore, consider matrix $\tbZ$ that includes all linearly independent columns in matrices $\bZ_K$ and $\hat{\bZ}_{\hat{K}}$,
so that 
$\Span\{\tbZ\} = \Span\{\bZ_K,\hat{\bZ}_{\hat{K}}\}$.
Let $\tK$ be the number of columns of matrix $\tbZ$. Then, one has
\begin{align*}
& \bPi_{\hat{\bZ},\hat{K}} \bPi_{\tbZ,\tK}=\bPi_{\tbZ,\tK}\bPi_{\hat{\bZ},\hat{K}} =\bPi_{\hat{\bZ},\hat{K}}, \\
& \bW_J = \bW_{J_1} + \bW_{J_3},\ \bW_{\hJ} = \bW_{J_1} + \bW_{J_2},\  
\bW_{\tJ} = \bW_{J_1} + \bW_{J_2} +\bW_{J_3}. 
\end{align*}
In order to obtain an upper bound for $\Del_2$ defined in \fr{Del1_2}, note that using notations above, we can rewrite $\Del_2$ as
\begin{align*}
\Del_2 
& =  {2\,\del}  \boeta^T(\bI_M \otimes \bS^T\bUp )[(\bPi_{\hat{\bZ},\hat{K}} \otimes \bW_{J_2})+(\bPi_{\hat{\bZ},\hat{K}} 
\otimes \bW_{J_1})-(\bPi_{\bZ,K} \otimes \bW_{J_1})-(\bPi_{\bZ,K} \otimes  \bW_{J_3})]\bg_* \\
& =   {2\,\del}  \boeta^T(\bI_M \otimes \bS^T\bUp )[(\bPi_{\hat{\bZ},\hat{K}}\otimes \bW_{J_2})+(\bPi_{\tbZ,\tK} \otimes  \bW_{J_1})+
(\bPi_{\bZ,K} \otimes  \bW_{J_3})] [(\bPi_{\hat{\bZ},\hat{K}} \otimes \bW_{J_2}) \\
& + (\bPi_{\hat{\bZ},\hat{K}} \otimes \bW_{J_1})-
(\bPi_{\bZ,K} \otimes \bW_{J_1})-(\bPi_{\bZ,K} \otimes  \bW_{J_3})]\bg_*\\
& =  {2\,\del}  \boeta^T(\bI_M \otimes \bS^T\bUp )[(\bPi_{\hat{\bZ},\hat{K}} \otimes \bW_{J_2}) + 
(\bPi_{\tbZ,\tK} \otimes  \bW_{J_1})+(\bPi_{\bZ,K} \otimes  \bW_{J_3})]
[(\bPi_{\hat{\bZ},\hat{K},\hat{J}} -\bPi_{\bZ,K,J})]\bg_*
\end{align*}
Using Cauchy inequality  and $2ab \leq 4 a^2 + b^2/4$, we obtain 
\begin{align} \label{del2_21_22}
& |\Del_2|\leq |\Del_{2,1}|+|\Del_{2,2}|, \quad  
|\Del_{2,2}|= 0.25\,  \|(\bPi_{\hat{\bZ},\hat{K},\hat{J}}\bg_*-\bPi_{\bZ,K,J}\bg_*)\|^2 \\
& |\Del_{2,1}|=4\, \del^2\, \|[(\bPi_{\hat{\bZ},\hat{K}} \otimes \bW_{J_2})+(\bPi_{\tbZ,\tK} \otimes  \bW_{J_1})
+(\bPi_{\bZ,K} \otimes  \bW_{J_3})](\bI_M \otimes \bUp\bS )\boeta\|^2 \nonumber
\end{align} 
Applying  Cauchy  Inequality to the term $\Del_{2,1}$ and using that  $J_2 \subseteq \hat{J}$ and $J_3 \subseteq J$  
we   rewrite 
%
\bes  
|\Del_{2,1}| \leq    12\del^2\, \lkv \|(\bPi_{\hat{\bZ},\hat{K}} \otimes (\bW_{\hat{J}}\bUp\bS ))\boeta\|^2
+\|(\bPi_{\tbZ,\tK} \otimes  (\bW_{J_1}\bUp\bS) )\boeta\|^2+\|(\bPi_{\bZ,K} \otimes  (\bW_{J}\bUp\bS) )\boeta]\|^2 \rkv
\ees
The upper bounds for the first and the third term in the inequality above can be obtained directly 
from Lemma~\ref{lem:lemma1}. For the second term, note that since $\tK \leq K + \hK$ and 
$J_1 \subseteq J$ and $J_1 \subseteq \hat{J}$ for any $\om \in \Om_{1\tau} \cap \Om_{2\tau}$ one has 
\begin{align} \label{eq3:lemma1}
\|(\bPi_{\tbZ,\tK} \otimes  (\bW_{J_1}\bUp\bS) )\boeta\|^2 & 
\leq   C_\psi^2\, \lkv 2K\,  \sum_{j \in J}\nu_j^2 + 2\hK\, \sum_{j \in \hat{J}} \nu_j ^2 \right.\\
\nonumber & \left. + 3 \, (\max_{j \in \hat{J}} \nu_j^2) \lfi  M \ln \hK  +|\hat{J}|\ln \lkr \frac{ne}{|\hat{J}|} \rkr +\ln (Mn) + 
 \tau \ln (\del^{-1})  \rfi \rkv
\end{align}
due to 
\bes
\tK \sum_{j \in J_1} \nu_j^2 \leq  K\,  \sum_{j \in J}\nu_j^2 +  \hK\, \sum_{j \in \hat{J}} \nu_j ^2.
\ees
Combining \fr{eq3:lemma1} with equations \fr{eq1:lemma1} and \fr{eq2:lemma1}, we obtain for any $\om \in \Om_{1\tau} \cap \Om_{2\tau}$ 
\begin{align}  \label{del21}
|\Del_{2,1}| & \leq  12 \del^2\,  C_\psi^2\, \lkv 4 \hK\, \sum_{j \in \hat{J}} \nu_j ^2 + 
4 K\,  \sum_{j \in J}\nu_j^2    + 3 (\max_{j \in J} \nu_j ^2)(\tau \ln n) \right. \\
\nonumber & \left. + 6\, (\max_{j \in \hat{J}} \nu_j ^2) \lfi M \ln \hK + |\hat{J}|\ln \lkr \frac{ne}{|\hat{J}|}\rkr  
+ \ln (M n) +  \tau \ln (\del^{-1})   \rfi  \rkv   
\end{align}


Now consider $|\Del_{2,2}|$ defined in \fr{del2_21_22}. Rewrite $|\Del_{2,2}|$ as 
$|\Del_{2,2}|= 0.25 \|(\bPi_{\hat{\bZ},\hat{K},\hat{J}}\bg_*-\bg_*)-(\bPi_{\bZ,K,J}\bg_*-\bg_*)\|^2$,
so that 
\bes
|\Del_{2,2}| \leq 0.5\,  \|(\bPi_{\hat{\bZ},\hat{K},\hat{J}}\bg_*-\bg_*)\|^2+ 0.5  \|(\bPi_{\bZ,K,J}\bg_*-\bg_*)\|^2.
\ees
Since $\hat{\bg}=\bPi_{\hat{\bZ},\hat{K},\hat{J}}\bGa\by$ and
\begin{align*}
 \|(\bPi_{\hat{\bZ},\hat{K},\hat{J}}\bGa\by-\bg_*)\|^2 &  =\|(\bPi_{\hat{\bZ},\hat{K},\hat{J}}(\bg_*+ \del \,\bGa \beps)-\bg_*)\|^2\\
 & =\|(\bI - \bPi_{\hat{\bZ},\hat{K},\hat{J}})\bg_*\|^2+ \del^2 \, \|\bPi_{\hat{\bZ},\hat{K},\hat{J}}\,  \bGa \beps  \|^2,
\end{align*}
we derive  
\be \label{Ineq_2}
\|\hat{\bg} -\bg_*\|^2  \geq \| \bPi_{\hat{\bZ},\hat{K},\hat{J}}\bg_*-\bg_* \|^2
\ee
Taking into account that   $\bg=\bPi_{\bZ,K,J}\bg_*$, so that $\|\bg-\bg_*\|^2=\| \bPi_{\bZ,K,J}\bg_*-\bg_* \|^2$,
we obtain 
\be \label{del22}
|\Del_{2,2}| \leq  0.5 \|\hat{\bg} -\bg_*\|^2+ 0.5  \|\bg-\bg_*\|^2.
\ee 
%
%
By combining upper bounds of $\Del_1$, $\Del_{2,1}$ and $\Del_{2,2}$,
we derive from \fr{eq:Del1} and  \fr{del21}-- \fr{del22} that 
for any $\om \in \Om_{1\tau} \cap \Om_{2\tau}$, an upper bound for $\Del$ can be written as 
\begin{align}   
\nonumber |\Del|  & \leq  0.5 \|\hat{\bg} -\bg_*\|^2+ 0.5  \|\bg-\bg_*\|^2 +   {2\, \del^2\, C_\psi^2} \lfi 
26   \hK\, \sum_{j \in \hat{J}} \nu_j ^2+ 24   K\,  \sum_{j \in J}\nu_j^2  \right.\\
\label{delta}
 & \left. + 39 \, (\max_{j \in \hat{J}} \nu_j ^2) \lkv M \ln \hK + |\hat{J}|\ln \lkr \frac{ne}{|\hat{J}|}\rkr  +\ln (M n) +  \tau \ln (\del^{-1})  \rkv  
 + 18(\max_{j \in J} \nu_j ^2)\tau \ln (\del^{-1}) \rfi  
\end{align}
Since it follows from \fr{vec_notations2} that 
$\|\widehat{\bG}-\bG_*\|_F ^2 = \|\hat{\bg}-\bg_*\|^2$, we obtain from \eqref{VLasso} that 
for any $\bG = \Pi_{\bZ,K,J} \bG_*$ on the set $\Om_{1\tau} \cap \Om_{2\tau}$ one has 
\begin{align}\label{Lasso_final}
\|\widehat{\bG} & -\bG_*\|_F ^2    \leq  3\|\bG-\bG_*\|_F ^2 + {2\,\del^2\, C_\psi^2}  \lfi  
48   K\,  \sum_{j \in J}\nu_j^2  + 36 (\max_{j \in J} \nu_j ^2) \tau  \ln \del^{-1}  + 52   \hK\, \sum_{j \in \hat{J}} \nu_j ^2   \right.\\
& \left.
 + 78  (\max_{j \in \hat{J}} \nu_j ^2) \lkv M \ln \hK  
  + |\hat{J}|\ln \lkr \frac{ne}{|\hat{J}|}\rkr  +\ln (M n) +  \tau \ln \del^{-1}  \rkv \rfi + 2 [\Pen(J,K) -  \Pen(\hat{J},\hat{K})] \nonumber 
\end{align}
%
%
Choose $\Pen(J,K)$ in the form \eqref{eq:penalty} and note that all terms containing $\hat{J}$ and $\hat{K}$  in \eqref{Lasso_final} cancel. 
Finally we obtained for any $\bG=\bW_J \bG_* \bPi_{\bZ,K}$ that with probability at least $1 - 2 \del^{\tau}$
\begin{align*} 
  \|\widehat{\bG}-\bG_*\|_F ^2  &\leq   3\|\bG-\bG_*\|_F ^2 +  {2\, \del^2\, C_\psi^2} \lfi  
48   K\,  \sum_{j \in J}\nu_j^2  + 36 (\max_{j \in J} \nu_j ^2) \tau  \ln n \rfi + 2\, \Pen(J,K) 
\end{align*}
which yields \eqref{eq:upper_bound}. 
\\


\subsection{Proof  of the  upper bounds for the error}

\noindent
{\bf Proof of Theorem~\ref{th:upper_bound_L}. } \ 
Since, when $j$ is growing, coefficients $\bTe_{jk}$ are decreasing while the values of $\nu_j$ 
are increasing according to \eqref{nu_j_cond},   the optimal set $J$ is of the form  $J = \lfi 1, \cdots, L \rfi$, 
so that $|J| = L$.  Then, we  find $(\hat{\bZ}, \widehat{\bG},\hat{L}, \hat{K})$ 
as a solution of optimization problem \eqref{opt_prob1} with the penalty given by expression \eqref{eq:penalty_new}.

Note that for the true number of classes $K_*$ with $N_k, k=1, \ldots, K_*$ elements in 
each class, $\bG$ are coefficients of each $f_m$ and $\bTe$ is the clustered version of 
those coefficients. It follows from \eqref{main_relation} that
\be \label{RMKn}
R (\hbf,\calS(r,\calA),M,K_*)  \leq  M^{-1} \|\widehat{\bG}-\bG_*\|_F ^2 + M^{-1} 
\sum_{k=1}^{K_*}\ N_k \sum_{j=n+1}^\infty \, \bTe_{jk}^2.
\ee 
Therefore, application of  the upper bound \eqref{eq:upper_bound_new} with a generic $L$,  $\bZ =\bZ_*$, $\bK = K_*$,
where $\bZ_*$ and $K_*$ are  respectively the true clustering matrix and the true number of classes, 
yields
\be \label{eq:oracle_here}
M^{-1}\, \|\widehat{\bG}-\bG_*\|_F ^2    \leq   3\, M^{-1}\,  \|\bW_J \bG_* \bPi_{\bZ_*,K_*} - \bG_*\|_F ^2 + 
4\, M^{-1}\,  \tPen(L,K_*) 
\ee
where $\tPen(L,K)$ is defined in \eqref{eq:penalty_new}.
Observe that 
\be \label{bias1}
\|\bW_J \bG_* \bPi_{\bZ_*,K_*}-\bG_*\|_F ^2 =  \|(\bW_J- \bI_n )\bG_*\|_F ^2=  \sum_{k=1}^{K^*} N_k   \sum_{j=L+1}^{n}\bTe_{jk}^2  
\ee 
where $N_k$ is the number of functions $f_m = h_k$ in the cluster $k$, $k=1, \cdots, K^*$, and $\bTe_{jk}$ are the true coefficients 
of those functions.   
Hence,  it follows from  \eqref{coef_cond} that 
\be \label{eq:tail_sum}
 \sum_{j=L+1}^{n}\bTe_{jk}^2 \leq \calA^2 L^{-2r}.
\ee
Since  $\displaystyle \sum_{k=1}^{K^*} N_k = M$, \eqref{bias1} and \eqref{eq:tail_sum} yield
\be \label{bias_L}
\|\bW_J \bG_* \bPi_{\bZ_*,K_*}-\bG_*\|_F ^2 \leq \calA^2 M L^{-2r}  
\ee
Moreover, it follows from \eqref{eq:tail_cond}    that
\bes
M^{-1}\ \sum_{k=1}^{K_*}\ N_k \sum_{j=n+1}^\infty \, \bTe_{jk}^2 \leq   \calA^2   n^{-2r} \asymp \del^2,
\ees
so that the last term in \eqref{RMKn} is   smaller than  $C\,  R (M,K_*,\del)$.

Now, consider the second term in  \eqref{eq:oracle_here}.  
Due to the condition \eqref{nu_j_cond}, one obtains
\bes 
\nu_L ^2   \leq \aleph_2 ^2\,  L^{2\ga} \exp\lkr  2\alpha L^\beta \rkr, \quad
\sum_{j=1}^L \nu_j ^2 \leq \aleph_2 ^2\,  L^{2\ga +1} \exp\lkr  2\alpha L^\beta \rkr. 
\ees
%
%
%
%
Denote 
\be \label{R1R2} 
R_1 \equiv R_1(K_*,\del) \asymp   K_*, \quad    R_2 \equiv R_2(M,K_*,\del) \asymp  M \ln K_*+\ln(\del^{-1}).   
\ee  
Therefore, it follows from  \eqref{eq:penalty}  and   
\eqref{eq:upper_bound} that, under condition \eqref{eq:M_n_rel}, 
\be  \label{up_bou1}
 \frac{\|\widehat{\bG}-\bG_*\|_F ^2}{M}   \leq \tilde{C}\,   \min_{L} \lfi   L^{-2r}  +  
 \frac{\del^2\, L^{2\ga}\, \exp (2\alpha L^\beta)}{M} \, \lkv  L R_1(K_*,\del) +   R_2(M,K_*,\del)  \rkv \rfi\ \ \ \ 
\ee
where $R_1(K_*,\del)$ and $R_2(M,K_*,\del)$ are defined in \eqref{R1R2} and 
$\tilde{C}$ depends only on $\mu$, $\calA$, $\aleph_2$, $C_\psi ^2$  and is independent of $M$, $L$, $\del$ and $K_*$.

In order to find the  minimum of the right hand side of \eqref{up_bou1}, denote
\be \label{opti_prob1} 
R(L,M,K_*,\del)=  L^{-2r} +  \del^2\, M^{-1}\, \exp\lkr  2\alpha L^\beta \rkr \lkv   L^{2\ga+1} R_1  + L^{2\ga} R_2 \rkv 
\ee
and observe that 
\be \label{upper}
M^{-1}\, \|\widehat{\bG}-\bG_*\|_F ^2   \leq \tilde{C}\,   \min_{L} R(L,M,K_*,\del) 
\ee 
where $L_{opt}$ is the value of $L$ minimizing the right hand side of \eqref{up_bou1}. Denote
\be \label{eq:Lopt1}
L_{1,opt} = \underset{L}{\operatorname{argmin}}\, [L^{-2r} +  \del^2\, M^{-1}\, \exp\lkr  2\alpha L^\beta \rkr     L^{2\ga+1} R_1],
\ee
\be \label{eq:Lopt2}
L_{2,opt} = \underset{L}{\operatorname{argmin}}\, [L^{-2r} +  \del^2\, M^{-1}\, \exp\lkr  2\alpha L^\beta \rkr     L^{2\ga} R_2].
\ee
It is easy to see that since the first terms in expressions \eqref{eq:Lopt1} and \eqref{eq:Lopt2}
are decreasing in $L$ while the second terms are increasing, the values $L_{1,opt}$ and $L_{2,opt}$ 
are such that those terms are equal to each other up to a multiplicative constant. 
Then, 
$R(L_{opt},M,K_*,\del) = \max \lfi L_{1,opt}^{-2r}, L_{2,opt}^{-2r} \rfi$, and, due to  
$\max(a,b) \asymp a+b$ for positive $a$ and $b$, we obtain
\be \label{eq:R_L_opt}
R(L_{opt},M,K_*,\del) \asymp L_{1,opt}^{-2r} + L_{2,opt}^{-2r}.
\ee 
Consider two cases.
\\

\noindent
{\bf Case 1:} $\al=\beta=0$.\ \ 
Direct calculations yield
\bes
L_{1,opt} \asymp \lkr M^{-1} \del^2 R_1\rkr ^ {-\frac{1}{2\ga +2r +1}}, \quad
L_{2,opt} \asymp \lkr M^{-1} \del^2 R_2\rkr ^ {-\frac{1}{2\ga +2r }},  
\ees
so that, due to \eqref{R1R2},  
\bes
L_{1,opt}   =  (M^{-1}\, \del^2  K_*)^{-\frac{1}{2\ga +2r +1}}, \quad 
L_{2,opt}   =  [\del^2 (\ln K_* + M^{-1} \ln \del^{-1})]^{-\frac{1}{2\ga +2r}} 
\ees
Then, by \eqref{eq:R_L_opt},  
\be \label{eq:RLopt_case1}
R(L_{opt},M,K_*,\del) \asymp (M^{-1}\, \del^2  K_*)^{\frac{2r}{2\ga +2r +1}}
+ [\del^2 (\ln K_* + M^{-1} \ln \del^{-1})]^{\frac{2r}{2\ga +2r}}.
\ee
Now, in order to obtain the  expression  \eqref{eq:RMKn_0}, note that if $K_* \geq 2$, then 
$\ln K_*$ dominates $M^{-1} \ln \del^{-1}$. If $K_* =1$, then \eqref{eq:RLopt_case1} 
can be re-written as 
$$
R(L_{opt},M,K_*,\del) \asymp \lkr\frac{\del^2}{M}\rkr^{\frac{2r}{2\ga +2r +1}}
\lkv 1 + \lkr\frac{\del^2}{M}\rkr^{\frac{2r}{(2\ga +2r +1)(2r + 2\ga)}} (\ln \del^{-1})^{\frac{2r}{2\ga +2r}} \rkv
\asymp \lkr\frac{\del^2  K_*}{M}\rkr^{\frac{2r}{2\ga +2r +1}}, 
$$
which yields \eqref{eq:RMKn_0}.
\\

\noindent
{\bf Case 2:} $\al>0, \beta>0$.\ \ 
Minimizing expressions in \eqref{eq:Lopt1} and \eqref{eq:Lopt2}, we obtain
\bes
L_{i,opt} \asymp  \lfi  \lkv \ln\lkr \frac{M }{\del^2 R_i}\rkr \rkv \rfi^{\frac{1}{\beta}},\quad i=1,2,
\ees
%
If $K_* \geq 2$, then $R_2 \geq R_1$. 
Taking into account that, under assumption \eqref{eq:M_n_rel}, for large $M$ and small  $\del$, 
$\ln \lkr  M  \del^{-2} \ln M  \rkr \asymp \ln \lkr M \del^{-2} \rkr $ and $\ln (Mn) \asymp \ln M$, we obtain  
\bes 
L_{1,opt} = \min \lfi  \lkv \ln\lkr \frac{M}{\del^2 K_*}\rkr \rkv ; 
\lkv \ln\lkr \frac{M}{\del^2 \ln M}\rkr \rkv\rfi^{\frac{1}{\beta}}
\asymp  \lkv \ln\lkr \frac{M}{\del^2 K_*}\rkr \rkv^{\frac{1}{\beta}}.
\ees 
Similarly, 
\bes
L_{2,opt} = \min \lfi  \lkv \ln\lkr \frac{1}{\del^2 \ln K_*}\rkr \rkv ; 
\lkv \ln\lkr \frac{M }{\del^2 \ln  M}\rkr \rkv\rfi^{\frac{1}{\beta}}
\asymp  \lkv \ln\lkr \frac{1}{\del^2 \ln K_*}\rkr \rkv^{\frac{1}{\beta}}, 
\ees
which, together with \eqref{upper} and \eqref{eq:R_L_opt}, yield  the  expression  \eqref{eq:RMKn_1}.
One can easily check that the case of $K_*=1$ leads to the same results.


\subsection{Proofs of the minimax lower bounds for the error}

\input{RP_Lowerbound_March22_2020}


\subsection{Proofs of the comparison of the risks with and without clustering}

\noindent
{\bf Proof of Corollary~\ref{cor:comparison}. } \ 
First observe that expressions \eqref{eq:noclust_upper_bound} are obtained 
directly from \eqref{eq:RMKn_0} and \eqref{eq:RMKn_1}  by setting $M=K_*=1$
since all functions belong to the same Sobolev ball \eqref{class_calB}.
In order to compare the upper bounds \eqref{eq:RMKn_0} and \eqref{eq:RMKn_1} 
obtained with clustering with the upper bound \eqref{eq:noclust_upper_bound}
derived without clustering, we consider several cases.

 \noindent
{\bf{Case 1} }  $\al=0$ , $\beta=0$.
\\
Expressions in \eqref{eq:comparison} are obtain by direct evaluation. Note that
the second expression in the case of $K_* \geq 2$ tends to zero as $M \to \infty$ 
since, due to \eqref{eq:M_n_rel}, $\ln K_* \leq \ln M \asymp \ln \del^{-1}$.
\\ 

\noindent
{\bf{Case 2} }  $\al>0$ , $\beta>0$.
\\
Note that, due to the condition \eqref{eq:M_n_rel}, 
\bes  
\ln (\del^{-2}) \leq \ln (M  \del^{-2} K_*^{-1}) \leq \ln M + \ln (\del^{-2}) \asymp  \ln (\del^{-2}),
\ees 
Also, for $K_* \geq 2$ and $\del^{-2} \geq e$, due to $\ln x \leq x/2$ for $x \geq 1$, obtain  
\bes
\ln \lkr \del^{-2}\, \ln (K_*^{-1}) \rkr = \ln (\del^{-2}) - \ln \ln K_* \geq \ln (\del^{-2}) - 0.5\,   \ln (\del^{-2}) \asymp   \ln (\del^{-2}),
\ees
which completes the proof.

\subsection{Proofs of supplementary statements}
\label{sec:suppl_proofs}


\noindent
{\bf Proof of Lemma  \ref{lem:lemma1}.\ }
Proof of Lemma~\ref{lem:lemma1}  is based on the following statement provided in  Gendre(2014)

\begin{lemma}\label{lem:lemma2} {\bf (Gendre (2014)).  } Let $\bA \in R^{p \times p}$ 
be a fixed matrix and $\beps \sim N(0,\bI_p)$. Then,  for any $x > 0$  one has
\be \label{eq:lemma2}
\PP \lfi \|\bA \beps\|^2 \geq \Tr (\bA^T \bA) +2\, \sqrt{\|\bA\|_{op} ^2 \Tr (\bA^T\bA) \, x}  +2\|\bA\|_{op} ^2 x \rfi \leq e^{-x}
\ee
\end{lemma}

\medskip

\noindent
Note that,  due to $2ab \leq a^2+b^2$, probability \eqref{eq:lemma2} can be re-written as 
\be\label{Gendre}
\PP(\|\bA\beps\|^2 \geq 2\|\bA\|_F ^2 +3\|\bA\|_{op} ^2 x) \leq e^{-x}
\ee
%
%
Consider $ \|[\bPi_{\bZ,K}  \otimes (\bW_{J} \Up \bS)]\,\boeta  \|^2$ with $\bZ,J,K$ fixed.
Note that, due to $\| \bPi_{\bZ,K} \|^2_{op}=1$, $\|\bS  \|^2_{op}\leq C^2_\psi$, $\| \bW_J \Up \|^2_{op}=\max_{j \in J} \nu_j ^2$ 
and $\| \bW_{J} \Up \|^2_F =  \sum_{j \in J}\nu_j ^2$, one has  
\be\label{opertor_norm_3}
\|( \bPi_{\bZ,K}  \otimes (\bW_{J} \Up \bS)) \|^2_{op} \leq \| \bPi_{\bZ,K} \|^2_{op} \| \bW_{J} \Up \|^2_{op}  \|\bS  \|^2_{op} 
\leq C^2_\psi \max_{j \in J} \nu_j ^2
\ee
\be\label{Frobenius_norm_4}
\|( \bPi_{\bZ,K}  \otimes (\bW_{J} \Up \bS))\boeta \|^2_F \leq \| \bPi_{\bZ,K} \|^2_{F} \| \bW_{J} \Up \|^2_F \|\bS  \|^2_{op} 
\leq K C^2_\psi \sum_{j \in J}\nu_j ^2
\ee
Now applying inequality  \eqref{Gendre} to 
$ \|[\bPi_{\bZ,K}  \otimes (\bW_{J} \Up \bS)]\, \boeta \|^2$ where $\boeta \sim N(0,\bI_{nM})$, we obtain for any $x>0$
\begin{align} 
& \PP \lfi \|( \bPi_{\bZ,K}  \otimes (\bW_{J} \Up \bS)) \boeta \|^2 \geq 2\|( \bPi_{\bZ,K}  \otimes (\bW_{J} \Up \bS)) \|^2_F 
+ 3\|( \bPi_{\bZ,K}  \otimes (\bW_{J} \Up \bS)) \|^2_{op}\, x \rfi \leq \nonumber \\
& \PP \lfi \|( \bPi_{\bZ,K}  \otimes (\bW_{J} \Up \bS)) \boeta \|^2 - C^2_\psi \lkv 2\, K  \sum_{j \in J}\nu_j ^2 
+ 3 x\,  \max_{j \in J} \nu_j^2  \rkv \geq 0  \rfi \leq e^{-x}. 
\label{main_ineq_lem1}
\end{align}
Setting $x=\tau \ln (\del^{-1})$  yields \eqref{eq1:lemma1}. 
Inequality \eqref{eq1:lemma1_new}  follows from \eqref{eq1:lemma1} since $\nu_j$ are growing with $j$ and 
$J = \{1, ..., L\}$.
\\

\noindent
In order to prove inequality \eqref{eq2:lemma1}, note that  for
\bes   
x (M,K,|J|,s) =M \ln K+|J|\ln (ne/|J|)+ \ln (M n) +  s,
\ees 
due to $\ln{n \choose j} \leq j\ln (\frac{ne}{j})$, one has 
\begin{align}
\sum_{\bZ,K,J} e^{-x (M,K,|J|,s)} & \equiv \sum_{K=1} ^M\, \sum_{j=1} ^n\, \sum_{|J| =j}\, \sum_{\bZ \in \calM(M,K)} e^{-x(M,K,j,s)} \nonumber \\
& =  \sum_{K=1}^M \sum_{j=1}^n {n \choose j} K^M e^{-x(M,K,j,s)}\nonumber \\ 
& \leq  \sum_{K=1}^M \sum_{j=1}^n  \lkr \frac{ne}{j}\rkr^j K^M e^{-x(M,K,j,s)} \leq  e^{-s} \label{eq:sum_ineq}
\end{align}
Therefore, by \eqref{main_ineq_lem1} and \eqref{eq:sum_ineq}, we obtain
\begin{align*}
&   \PP \lkr \|(\bPi_{\hat{\bZ},\hat{K}}  \otimes (\bW_{\hat{J}} \Up \bS)) \boeta \|^2  - 
2\|( \bPi_{\hat{\bZ},\hat{K}}  \otimes (\bW_{\hat{J}} \Up \bS)) \|^2_F - 
3\|( \bPi_{\hat{\bZ},\hat{K}}  \otimes (\bW_{\hat{J}} \Up \bS)) \|^2_{op}\, x(M,\hat{K},|\hat{J}|,s)  \geq 0 \rkr \leq   \\
& \sum_{\bZ,K,J} \PP \lkr \|( \bPi_{\bZ,K}  \otimes (\bW_{J} \Up \bS)) \boeta \|^2 -  C^2_\psi \lkv 2K  \sum_{j \in J}\nu_j ^2 + 
3\, x(M,K,|J|,s) \lkr \max_{j \in J} \nu_j^2 \rkr   \rkv \geq 0  \rkr \leq \\
& \sum_{\bZ,K,J} e^{-x (M,K,|J|,s)}  \leq  e^{-s}. 
\end{align*} 
Setting $s=\tau \ln (\del^{-1})$  yields \eqref{eq2:lemma1}.

Similarly, in order to prove \eqref{eq2:lemma1_new}, choose $J = \{1, ..., L\}$,
$x (M,K,|J|,s) = M \ln K +  \ln (M n) +  s$,  and replace \eqref{eq:sum_ineq} by
\begin{align*}
\sum_{\bZ,K,J} e^{-x (M,K,|J|,s)} & \equiv \sum_{K=1} ^M\, \sum_{L=1} ^n\,   \sum_{\bZ \in \calM(M,K)} e^{-x(M,K,L,s)}   \\
& \leq  \sum_{K=1}^M  n\,  K^M e^{-x(M,K,L,s)}  \leq  e^{-s}  
\end{align*}


\medskip


\noindent
{\bf Proof of Lemma  \ref{lem:fd}.\ }
By using \eqref{d_cond}, $K \geq 2$ and $0< d \leq 1/9$
\begin{align*}
\ln K - 4d \ln (K e/d) & = \ln K - 4[d \ln (K)+d-d\ln d]\\
                       & \geq  \ln K-4d\ln K-\frac{4}{9}\ln 2 \\
                       & \geq \frac{5}{9}\ln K - \frac{4}{9}\ln K \geq \frac{\ln K}{9}.  
\end{align*}


%% file: RP_Lowerbound_March22_2020.tex

\noindent
{\bf Proof of Theorem~\ref{th:lower_bound_L}. } \  
Since the estimation error is comprised of the error due to nonparametric estimation and to clustering, we 
consider two cases here. 
\\

\noindent
\underline{\bf Lower bound for the error due to clustering}.  
\\
Let $K \geq 2$ be the fixed number of classes.  
Consider a subset $\calZ (M,K) \subset \calM (M,K)$ of the set of all clustering matrices which contains all matrices that cluster 
$\frac{M}{K}$ vectors into each class. By Lemma 5 in Pensky (2019) 
with $\ga =1$, obtain that the cardinality of the set $\calZ (M,K)$ is
\be \label{low_EQU1}
|\calZ(M,K)|= M!\Big/[(M/K)!]^K  \geq \exp\lkr M \ln K/4 \rkr
\ee
Let set $J$  be of the form  $J= \{L_1,...,L_2\}$ where $1 \leq L_1< L_2 \leq n$ and $n = [\del^{-2}]$.
Choose   $\bTe_{jk} = 0$ if $j \notin J$. 
In what follows, we use the Packing Lemma (Lemma~4 of Pensky (2019)):

\begin{lemma} \label{lem:packing} {(The Packing lemma).}
Let $\calZ (M,K) \subseteq \calM (M,K)$ be a collection of clustering matrices
and   $q$ be a positive constant. 
Then, there exists a subset $\calS_{M,K} (q) \subset \calZ (M,K)$ such  that 
for $\bZ_1, \bZ_2 \in \calS_{M,K} (q)$ one has 
$\| \bZ_1 - \bZ_2 \|_H =  \| \bZ_1 - \bZ_2 \|^2 _F \geq q$
and 
$\ln |\calS_{M,K} (q)| \geq  \ln|\calZ (M,K)| - q \ln (M  K e/q)$.
\end{lemma}

\noindent
Apply this lemma with $q= d M$, $0<d<1/4$.
Then, by \eqref{low_EQU1}, derive 
\bes
\ln|\calS_{M,K}\lkr  d M \rkr| \geq   M\,  \lkv \ln K - 4 d   \ln (K e/d) \rkv\big/ 4. 
\ees
Use the following statement:

\begin{lemma} \label{lem:fd}
If $K \geq 2$ and $d$ is such that 
\be \label{d_cond}
d - d \ln d \leq (\ln 2)/9, \quad d \leq 1/9,
\ee 
then $\ln K - 4d \ln (K e/d) \geq (\ln K)/9$.
\end{lemma}

\noindent
It is easy to calculate that, e.g.,   $d = 0.0147$ satisfies the condition \eqref{d_cond}. 
Then, for $d$ obeying  \eqref{d_cond}, one has
\be \label{eq:clust_set}
\ln|\calS_{M,K} (dM)| \geq \frac{M}{36 } \ln K, \quad \|\bZ_1-\bZ_2\|_H \geq dM \  \mbox{for any}
\  \bZ_1, \bZ_2 \in \calS_{M,K} (dM),\ \bZ_1 \neq \bZ_2
\ee 
Consider a collection of binary vectors ${\bom} \in \{0,1\}^{|J|} $. 
By  Varshamov-Gilbert bound lemma, there exists a subset $\calW$ of those vectors such that, for any 
${\bom},{\bom'} \in \calW$  such that ${\bom} \neq {\bom'}$ one has
$\|{\bom}-{\bom'}\|_H \geq  |J|/8$ and  $\ln {|\calW|} \geq |J| \ln (2)/8$.
Choose a subset $\calW_K$ of $\calW$ such that $|\calW_K|=K$. This is possible if 
$K \leq 2^{|J|/8}$ which is equivalent to $|J| \ge 8\,\ln K/\ln 2$.
Consider a set of vectors $\bw \in \{0,1\}^{n}$ obtained by  packing ${\bom}$ with zeros for components not in $J$. Then 
\be \label{set_G_K}
\calW_K =\lfi {\bw_1},...,{\bw_K}\in \{0,1\}^{n}: \, \|\bw_i\|_0 \leq |J|,\  \|\bw_i-\bw_j\|_0 \geq  |J|/8,\ i \neq j \rfi
\ee 
Define matrix $\bW$ with columns $\bw_k$,  $k=1,...,K$.
%
Finally, form the set $\calG_{M,K}$ of matrices $\bG$ of the form
\bes
\calG_{M,K}=\left\{\bG \in R^{n \times M}: \bG=\te \bW \bZ^T , \bZ \in  \calS_{M,K} (dM) \right\}
\ees
where $d$ satisfies \eqref{d_cond} and  $\te >0$   depends on $M$,$\del$ and $K$. 
Note that, due to \eqref{eq:clust_set}, one has
\be \label{card_valG}
\ln |\calG_{M,K}| \geq (M \ln K)/36
\ee  
Let $\bZ_1,\bZ_2 \in \calS_{M,K}$ be two clustering matrices. 
Set  $\bG_1=  \te \bW \bZ_1^T$  $\bG_2=  \te \bW \bZ_2^T$, so that $\bG_1,\bG_2 \in \calG_{M,K}$.
Since for any $i,i'$ one has $\|{\bw}_i-{\bw}_{i'}\|_0 =\|{\bw}_i-{\bw}_{i'}\|^2$, derive that
\begin{align} 
& \|\te \bW \lkr\bZ_1-\bZ_2\rkr^T\|_F^2
 =\sum_{m=1}^{M} \sum_{j=1}^{n}  \te^2 \lkv \lkr\bw_{z_1 {\lkr m\rkr}}\rkr_j- \lkr\bw_{z_2 {\lkr m\rkr}}\rkr_j\rkv^2 = \nonumber\\
&=\te^2 \sum_{m=1}^{M} \|{\bw}_{z_1 \lkr m\rkr}-{\bw}_{z_2 \lkr m\rkr}\|^2 \geq \#\{m:z_1\lkr m\rkr \neq z_2\lkr m\rkr\}\,  \te^2 |J|/8. 
\label{lowb1} 
\end{align} 
On the other hand, observe that for $\bZ_1,\bZ_2 \in \calS_{M,K}$  one has 
\bes
 \#\{m:z_1\lkr m\rkr \neq z_2\lkr m\rkr\} = 0.5\, \|\bZ_1-\bZ_2\|_H  \geq d M/2.
\ees
Therefore, the last two inequalities yield for any $\bG_1,\bG_2 \in \calG_{M,K}$
\be \label{low_eq2} 
\|\bG_1 - \bG_2\|_F^2 \geq  d \, \te^2 |J| M/16. 
\ee 
Now, it is easy to calculate that for any $\bG_1,\bG_2 \in \calG_{M,K}$ and the corresponding probability measures 
$P_{\bG_1}$ and $P_{\bG_2}$ associated with $\bY = \bUp^{-1} \bG_i + \del \bE$, $i=1,2$, in \eqref{main_seq}, 
one has the following inequality for the Kullback-Leibler divergence between 
$P_{\bG_1}$ and $P_{\bG_2}$:
\be \label{eq:Kullback}
K\lkr P_{\bG_1},P_{\bG_2} \rkr \leq \frac{1}{2\del^2 C_\psi^2} \| \bUp^{-1} \lkr \bG_2-\bG_1\rkr  \|_F ^2
\ee 
Since $\bG_1= \te \bW \bZ_1$, $\bG_2 = \te \bW \bZ_2$,  we obtain
\be  \label{low_eq3}
 \| \bUp^{-1} \lkr \bG_2-\bG_1\rkr  \|_F ^2 \leq \te^2 \, \| \bZ_2-\bZ_1 \|_{op} ^2 \, \| \bUp^{-1} \bW \|_F ^2
\ee 
Note that $\calS_{M,K} (dM) \subset \calZ (M,K)$, so that  for any $\bZ \in   \calS_{M,K} (dM)$ one has  $\bZ^T\bZ= (M/K)\,\bI_K$,
hence $\|\bZ\|_{op}= \sqrt{M/K}$.  Then,  $\|\bZ_1-\bZ_2\|_{op}^2 \leq 4M/K$. Also, 
due to $J= \{L_1,...,L_2\}$ and condition \eqref{nu_j_cond}, one has
\be \label{sum_nuj}
\sum_{j \in J}  \nu_j^{-2}  \leq  \aleph_1^{-2} |J|\, L_1 ^{-2\ga} \exp\lkr  -2 \alpha  L_1 ^\beta  \rkr.
\ee
Since 
$\| \bUp^{-1} \bW \|_F ^2 \leq  \sum_{k=1}^{K} \sum_{j \in J}  \nu_j ^{-2}$, obtain
\be \label{low_eq4}
K\lkr P_{\bG_1},P_{\bG_2} \rkr \leq \frac{2}{\del^2 \aleph_1^{2} C_\psi^2}\ \te^2   |J| M\,   L_1 ^{-2\ga} \exp\lkr  -2 \alpha  L_1 ^\beta  \rkr.
\ee
Finally, due to condition \eqref{coef_cond}, one needs 
$\te^2 \sum_{j\in J}   (j+1)^{2r} \leq  \calA^{2}$, 
so that we can choose 
\be \label{low_eq9}
\te^2 = \calA^2 |J|^{-1} L_2 ^{-2r}
\ee
In order to apply  Theorem 2.5  of Tsybakov (2009) with $\al=1/9$, we need 
$K\lkr P_{\bG_1},P_{\bG_2} \rkr \leq  \ln |\calG_{M,K}|/9$ which, 
due to \eqref{eq:clust_set}, is guaranteed by 
\be \label{Kullback_cond}
\frac{\te^2  |J|}{ \del^2 \aleph_1^{2} C_\psi^2}   L_1 ^{-2\ga} \exp\lkr  -2 \alpha  L_1^\beta  \rkr \leq  \frac{\ln K}{648}. 
\ee
If inequality  \eqref{Kullback_cond} holds, then application of Theorem 2.5  of Tsybakov (2009)    yields
that, with  probability at least 0.1, one has \eqref{eq:low_bou_main}
 where, due to \eqref{eq:fun_est_error} and  \eqref{low_eq2},
\be \label{lower_bound}
R_{\min} (M,K_*,\del) = \te^2 |J|.
\ee
Consider   $L_1=L/2+1$ and $L_2=L$, so that  
\be \label{theta_values}
\te^2 \asymp  L^{-(2r+1)}, \quad R_{\min} (M,K_*,\del) \asymp L^{-2r}.
\ee
If $\al=0$ , $\beta=0$, then,  by \eqref{theta_values},  inequality \eqref{Kullback_cond} holds if 
$L \asymp  \lkr  \del^2  \, \ln K \rkr ^{-\frac{1}{2r+2\ga}}.$
Hence,
\be \label{lower_1a}
R_{\min} (M,K_*,\del) \gtrsim \lkr  \del^2\, \ln K_* \rkr^{\frac{2r}{2r+2\ga}}.
\ee
If  $\al>0$, $\beta>0$,  then inequality  \eqref{Kullback_cond} holds if
$L ^{-(2\ga+2r)} \exp\lkr  -2 \alpha  L ^\beta  \rkr \lesssim   \del^2 \ln K$, so that
$L \asymp \lkv \ln \lkr \frac{1}{ \del^2 \ln K} \rkr \rkv ^{\frac{1}{\beta}}$. Therefore,   
\be \label{lower_2a}
R_{\min} (M,K_*,n) \gtrsim \lkv \ln \lkr \frac{1}{ \del^2 \ln  K_*} \rkr \rkv ^{-\frac{2r}{\beta}}.
\ee
 \\

  
\noindent
\underline{\bf Lower bound for the error due to estimation}.  
\\
Let, as before, $n= [\del^{-2}]$ and $J= \{L_1,...,L_2\}$ where $1 \leq L_1< L_2 \leq n$.
Consider a set of binary vectors $\bom \in \{0,1\}^{|J| K}$ and set $N = |J| K$. 
Complete  vectors $\bom$ with zeros to obtain vectors $\bw \in \{0,1\}^{nK}$. 
By  Varshamov-Gilbert  lemma, there exists a subset $\calB$ of those vectors such that for any $\bw,\bw'\in \calB$ such that $\bw\neq \bw'$
one has  $\|\bw-\bw'\|_H \geq  N/8$ and $\ln|\calB| \geq  N \ln (2)/8$.  
Pack vectors $\bw$ into matrices $\bW \in \{0,1\}^{n \times K}$. 
Denote the set of those matrices by $\calW$ and observe that 
\be \label{Varsh_Gil}
\|\bW_1 -\bW_2\|_F^2 \geq  N/8 \quad  \mbox{for all} \quad \bW_1, \bW_2 \in \calW,\  \bW_1 \neq \bW_2;  \qquad   \ln|\calW| \geq  (N \, \ln 2)/8.
\ee 
Let $\bZ$ be the clustering matrix that corresponds to uniform sequential clustering, 
$M/K$ vectors per class.  Finally, form the set $\calG_{M,K}$ of matrices $\bG$ of the form
\bes
\calG_{M,K}=\left\{\bG \in R^{M \times K}: \bG= \te \bW \bZ^T, \quad  \bW \in \calW   \right\}
\ees
where $\te >0$   depends on $M$,$\del$ and $K$.
Then, for any $\bG_1, \bG_2 \in \calG_{M,K}$, $\bG_1 \neq \bG_2$,  due to $\bZ^T\bZ=(M/K)\, \bI_K$ and 
\eqref{Varsh_Gil}, obtain 
\be \label{low_eq16}
\|(\bG_1-\bG_2) \|^2_F   =\te^2 \|(\bW_1-\bW_2)\bZ^T\|_F^2  
  = \frac{\te^2 M}{K}\,  \|\bW_1-\bW_2\|^2_F    \geq \frac{\te^2 M N}{8K}  
\ee 
Now, since  $\bG_1= \te \bW_1 \bZ$ and  $\bG_2 = \te \bW_2 \bZ$, using formula \eqref{eq:Kullback}, derive that 
\bes 
K \lkr P_{\bG_1},P_{\bG_2} \rkr \leq \frac{  \te^2}{2\del^2 C_\psi^2} \| \bUp^{-1} \lkr \bW_2-\bW_1\rkr \|_F ^2\, \|\bZ \|_{op} ^2 
\ees
Recalling that $\|\bZ \|_{op} ^2= M/K$ and $\| \bUp^{-1} \lkr \bW_2-\bW_1\rkr \|_F ^2 \leq \sum_{k=1} ^K \sum_{j \in J} \nu_j ^{-2}$,
and using \eqref{sum_nuj}, arrive at 
\bes
K \lkr P_{\bG_1},P_{\bG_2} \rkr \leq \frac{ M \te^2}{2\del^2\, \aleph_1^2\, C_\psi^2 }\, |J|\, L_1 ^{-2\ga} \exp\lkr  -2 \alpha  L_1 ^\beta  \rkr.
\ees
In order to apply  Theorem 2.5  of Tsybakov (2009)with  $\al=1/9$, we need 
$K\lkr P_{\bG_1},P_{\bG_2} \rkr \leq (1/9) \ln |\calG_{M,K}|$ which, 
due to \eqref{Varsh_Gil}, is guaranteed by 
\be \label{Kullback_cond_new}
\frac{\te^2  M}{ \del ^2 \aleph_1^{2} C_\psi^2}   L_1 ^{-2\ga} \exp\lkr  -2 \alpha  L_1^\beta  \rkr \leq  \frac{K}{36}. 
\ee
If inequality  \eqref{Kullback_cond_new} holds, then application of Theorem 2.5  of Tsybakov (2009)   yields
that, with  probability at least 0.1, one has \eqref{eq:low_bou_main},
where, due to \eqref{eq:fun_est_error} and  \eqref{low_eq16}, 
\be \label{lower_bound_est}
R_{\min} (M,K_*,\del) \gtrsim  \te^2   |J|
\ee
Now, as before, we consider two choices of $L_1$ and $L_2$: $L_1= L_2=L$ and $L_1=L/2+1$ , $L_2=L$ leading to the   
values of $\te^2$ given by \eqref{theta_values}. 
Again, we consider the cases of $\al= \beta=0$ and $\al>0$, $\beta>0$ separately. 
\\

\noindent
{\bf{Case 1:} }  $\al=0$ , $\beta=0$, $L_1=L/2+1$ , $L_2=L$, $|J|=L/2$.
\\
Since $L_1 \asymp L_2 \asymp |J| \asymp L$,  inequality \eqref{Kullback_cond_new} holds if
$L \asymp  \lkr  \del^2  M^{-1}\, K \rkr ^{-\frac{1}{2r+2\ga+1}} $
and 
\be \label{lower_3b}
R_{\min} (M,K_*,\del) \gtrsim   \lkr  \del^2\, M^{-1}\,  K  \rkr^{\frac{2r}{2r+2\ga+1}}.
\ee

\noindent
{\bf{Case 2:} }  $\al>0$, $\beta>0$, $L_1= L_2=L$,   $|J|=1$.
\\
Plugging the first expression from \eqref{theta_values} into \eqref{Kullback_cond_new}, derive that 
$L ^{-(2\ga+2r)} \exp\lkr  -2 \alpha  L ^\beta  \rkr \lesssim   \del^2  M^{-1} K$, so that
$L \asymp \lkv \ln \lkr \frac{M}{\del^2  K} \rkr \rkv ^{\frac{1}{\beta}}$. Therefore,   
\be \label{lower_4a}
R_{\min} (M,K_*,\del) \gtrsim \lkv \ln \lkr \frac{M}{ \del^2  K} \rkr \rkv ^{-\frac{2r}{\beta}}
\ee
 \\

\noindent
 Now, in order to obtain the expressions for the lower bounds, we find the maximum of 
\eqref{lower_1a}  and \eqref{lower_3b} if  $\al=0$ , $\beta=0$,
and of \eqref{lower_2a} and  \eqref{lower_4a}  if  $\al>0$ , $\beta>0$.


%% file: RP_March22_2020.bbl
\begin{thebibliography}{10}

\bibitem{abram}
Abramovich, F. and Silverman, B. W. (1998).  
Wavelet decomposition approaches to statistical inverse problems. 
{\it Biometrika}, {\bf 85}, 115--129.

\bibitem{abr_pen}
Abramovich, F., De Canditiis, D. and Pensky, M. (2018). 
Solution of linear ill-posed problems by model selection and aggregation.
{\it Electronic Journal of Statistics}, {\bf 12},   1822--1841.


\bibitem{alquier}
Alquier, P., Gautier, E. and Stoltz, G. (2011). 
 {\it  Inverse Problems and High-Dimensional Estimation}, Springer-Verlag, Berlin.


\bibitem{arnold}
Arnold, A., Reichling, S., Bruhns, O. T., and Mosler, J. (2010).
Efficient computation of the elastography inverse problem by combining variational 
mesh adaption and a clustering technique.
{\it Phys Med Biol.}, {\bf 55}, 2035-2056. 



\bibitem{bezdek}
Bezdek, J. C. and  Pal, S. K. (1992). 
{\it Fuzzy models for pattern recognition methods that search for 
structures in data}, IEEE Press, New York.


\bibitem{bissantz}
Bissantz, N., Hohage, T., Munk, A. and Ruymgaart, F. (2007).  
Convergence rates of general regularization methods for statistical 
inverse problems and applications.
{\it SIAM J. Numer. Anal.}, {\bf 45}, 2610-2636.


\bibitem{blanchard}
Blanchard, G.,  Hoffmann, M. and Reis, M. (2018).
Early stopping for statistical inverse problems via truncated SVD estimation.
{\it Electron. J. Statist.}, {\bf 12}, 3204-3231.

\bibitem{cohen}
Cohen, A., Hoffmann, M. and Reis, M. (2004).
Adaptive wavelet Galerkin methods for linear inverse problems.  
{\it SIAM Journ. Numer. Anal.}, {\bf 42}, 1479-1501.


\bibitem{com1}
Comte, F., Cuenod, C. A., Pensky, M. and   Rozenholc, Y. (2017).
Laplace deconvolution on the basis of   time domain data  
and its application to Dynamic Contrast Enhanced imaging.
{\it Journ. Royal Stat. Soc., Ser.B.}, {\bf 79},   69-94.

\bibitem{deng}
Deng, Z., Chung, F. L. and Wang, S. (2011).
Clustering-Inverse: A Generalized Model for 
Pattern-Based Time Series Segmentation. 
{\it Journal of Intelligent Learning Systems and Applications},  {\bf 3}, 26-36. 


\bibitem{donoho}
Donoho, D. L. (1995). Nonlinear solution of linear inverse problems
by wavelet-vaguelette decomposition {\it Applied and Computational
Harmonic Analysis}, {\bf 2}, 101--126.



\bibitem{donjohn}
Donoho,  D. L.  and Johnstone,  I. M.   (1994).  Ideal spatial adaptation by wavelet shrinkage.  
{\it Biometrika }{\bf 81}, 425--456. 


\bibitem{engl}
  Engl, H. W., Hanke, M. and  Neubauer, A. (2000).
{\it Regularization of Inverse Problems}, Kluwer Academic Publishers, Netherlands. 



\bibitem{fraix}
Fraix-Burnet, D. and  Girard, S.  (2016).
{\it Statistics for Astrophysics Clustering and Classification}, EDP Sciences.

\bibitem{gendre}
Gendre, X. (2014)
Model selection and estimation of a component in additive regression.
{\it ESAIM: Probability and Statistics}, {\bf 18}, 77--116.

 
\bibitem{gupta}
Gupta,  A. K. and  Nagar, D. K. (1999).
{\it Matrix Variate Distributions},
Chapman \& Hall/CRC, Boca Raton.




 


 
\bibitem{klopp}
Klopp, O.,  Lu Y.,   Tsybakov, A. B. and Zhou, H. H. (2019). 
Structured matrix estimation and completion.
{\it  Bernoulli}, {\bf   25},   3883--3911.



\bibitem{kurum}
K\"{u}r\"{u}m, E., Weber, G. W. and Iyigun, C. (2018).
Early warning on stock market bubbles via methods
of optimization, clustering and inverse problems.
{\it  Annals of Operations Research}, {\bf 260}, 293-320.


\bibitem{mallat} 
 Mallat, S.  (2009).
{\it A Wavelet Tour of Signal Processing. The Sparse Way.}
3rd Edition. Academic Press, New York.
 




\bibitem{pen_lasso}
Pensky, M. (2016).
Solution of linear ill-posed problems using overcomplete dictionaries.
{\it  Annals of Statistics}, {\bf 44},    1739--1764.


\bibitem{pen}
Pensky, M. (2019). Dynamic network models and graphon estimation.  
{\it  Annals of Statistics}, {\bf 47},  2378--2403.


\bibitem{starck}
Starck, J. L. and Pantin, E. (2002).
Deconvolution in Astronomy : A Review.
{\it Publ.  Astronom. Soc. of the Pacific}, {\bf 114},  1051-1069.



 \bibitem{tsybakov}
Tsybakov, A. B. (2009).
{\it Introduction to Nonparametric Estimation}, Springer, New York. 

 
 
\end{thebibliography}
